\documentclass[10pt, amsthm, reqno, tbtags]{elsartcld}
\usepackage{mathrsfs}
\usepackage{amsmath}
\usepackage{epsfig}
\usepackage{amsfonts}
\usepackage{amssymb}
\usepackage{amscd}
\usepackage[all]{xy}
\usepackage{rotating}
\usepackage{lscape}
\usepackage{fancyhdr}
\usepackage{amsbsy}
\usepackage{verbatim}
\usepackage{moreverb}

\newcommand{\ra}{\rightarrow}

\newcommand{\xra}{\xrightarrow}
\newcommand{\xla}{\xleftarrow}

\newcommand{\D}{{\rm D}}
\newcommand{\R}{{\mathbb R}}

\newcommand{\sm}{\wedge}

\newcommand{\Q}{{\mathbb{Q}}}
\newcommand{\id}{{\rm id}}

\newcommand{\tensor}{\otimes}

\newcommand{\Spin}{{Spin}}
\newcommand{\Spinc}{{Spin^c}}

\newcommand{\Z}{{\mathbb Z}}
\newcommand{\g}{{\mathfrak g}}
\newcommand{\ttt}{{\mathfrak t}}
\newcommand{\CP}{{\mathbb{C} \rm{P}}}
\newcommand{\C}{{\mathbb{C}}}
\newcommand{\G}{{\mathbb G}}
\newcommand{\OO}{{\mathcal{O}}}
\newcommand{\lb}{\langle}
\newcommand{\rb}{\rangle}
\newcommand{\HH}{{\mathcal{H}}}
\newcommand{\GL}{GL}

\newcommand{\bdot}{{\displaystyle{\cdot}}}

\newcommand{\nn}{\nonumber}
\newcommand{\nid}{\noindent}

\newcommand{\sss}{\scriptscriptstyle}

\renewcommand{\H}{{\rm H}}

\theoremstyle{theorem}
\newtheorem{theorem}{Theorem}[section]
\newtheorem{propos}[theorem]{Proposition}
\newtheorem{lemma}[theorem]{Lemma}

\theoremstyle{definition}

\theoremstyle{remark}
\newtheorem*{conject}{Conjecture}

 \DeclareMathOperator{\ind}{ind}
 
\DeclareMathOperator{\Tor}{Tor}
\DeclareMathOperator{\Ext}{Ext}
\DeclareMathOperator{\colim}{colim}
\DeclareMathOperator{\Aut}{Aut}

\DeclareMathOperator{\TMF}{TMF}

\DeclareMathOperator{\Thom}{Th}
\DeclareMathOperator{\Sym}{Sym}

\CompileMatrices

\begin{document}

\begin{frontmatter}

\title{\hspace*{.16in}{\mbox{On the Twisted $K$-Homology of Simple Lie Groups}}
\vspace{-36pt}
}
\author{\hspace*{.75in}Christopher L. Douglas\thanksref{cld}}
\address{\hspace*{.75in}
Department of Mathematics, M. I. T., Cambridge, MA
02139, U.S.A.}
\thanks[cld]{The author was supported by a National Science
Foundation Graduate Research Fellowship}

\begin{abstract}
\setlength{\textwidth}{6.3in}
\begin{minipage}{\textwidth}
We prove that the twisted $K$-homology of a simply connected simple Lie
group $G$ of rank $n$ is an exterior algebra on $n-1$ generators
tensor a cyclic group.  We give a detailed description of the order of
this cyclic group in terms of the dimensions of irreducible representations
of $G$ and show that the congruences determining this cyclic order lift
along the twisted index map to relations in the twisted $\Spinc$ bordism
group of $G$.
\end{minipage}

\end{abstract}

\end{frontmatter}

\setcounter{tocdepth}{3}
\tableofcontents


\section{Introduction}

By way of motivation we present six interpretations of twisted $K$-theory.
These interpretations inform the methods and perspectives adopted in
the paper but are otherwise unnecessary for what follows. We then summarize 
our results on the twisted $K$-homology of simple Lie groups and overview our 
main techniques, namely the twisted Rothenberg-Steenrod spectral sequence,
Tate resolutions, Bott generating varieties, and twisted $\Spinc$ bordism.

\vspace{-10pt}

\subsection{Six Interpretations of Twisted $K$-Theory}

\vspace{-5pt}

\subsubsection{1-Dimensional Elements in Elliptic Cohomology}

A twisting on a space $X$ of a cohomology theory represented by
a spectrum $R$ is a bundle of spectra on $X$ with fibre $R$ and
the associated twisted cohomology of $X$ is given by the homotopy 
classes of sections of this bundle.  Such twistings are classified
by maps from $X$ to the classifying space $B\Aut R$ of homotopy
automorphisms of the spectrum $R$.  If $R$ is an $A_{\infty}$
ring spectrum, the classifying space $B\GL_1 R$ of homotopy units
in $R$ maps to $B\Aut R$ and thereby classifies a subset of the
twistings---we refer to these twistings as elementary.

The classifying space $B\GL_1 \H\C$ for elementary
twistings of ordinary cohomology with complex coefficients is $B\C^*$;
here $\H\C$ denotes the Eilenberg-MacLane spectrum for $\C$, and
both $\C$ and $\C^*$ have the discrete topoology.
There is a map $B\C^* \ra \Z \times BU$ of this
classifying space into the representing space for $K$-theory;
any twisting $X \ra B\GL_1 \H\C$ for the ordinary cohomology of $X$
therefore determines a $K$-theory class on $X$.  Of course, there is a natural
geometric interpretation of the $K$-theory classes arising in this way,
namely as the classes represented by flat line bundles on $X$.  The
twisted cohomology of $X$ is simply the cohomology of $X$ with coefficients 
in the line bundle, reinterpreted as the homotopy classes of sections of
an associated $\H\C$ bundle.

The classifying space $B\GL_1 K$ for elementary twistings of complex $K$-theory
splits, as an infinite loop space, as $T \times S$.  The factor $T$ is a
$K(\Z,3)$ bundle over $K(\Z/2,1)$ which splits as a space but has nontrivial
infinite loop structure classified by $\beta Sq^2 \in H^3(\H(\Z/2);\Z)$.
There is a natural infinite loop map $T \ra \TMF$ from
$T$ to the representing space for topological modular forms, and
so by projecting through $T$ a map $B\GL_1 K \ra \TMF$.  In particular
an elementary twisting of $K$-theory for $X$ determines a $\TMF$-class on $X$.
(Notice that $\TMF$ is the analog of real $K$-theory, that is of $KO$, and so 
the map $B\GL_1 K \ra \TMF$ corresponds to the composite 
$B\GL_1 \H\C \ra \Z \times BU \ra \Z \times BO$; it is not known whether
there exists an appropriate factorization $B\GL_1 K \ra E \ra \TMF$ for every
elliptic spectrum $E$.)  The geometric interpretation of these
$\TMF$ classes is simplified if we restrict our attention to those
classes coming from twistings involving only the $K(\Z,3)$ factor of $T$.
Such a twisting is determined by a map $X \ra K(\Z,3)$ or equivalently
by a $BS^1$ bundle on $X$.  We think of this bundle as a stack locally
isomorphic to the sheaf of line bundles on $X$ and as such as a
1-dimensional 2-vector bundle on $X$.  In this sense we imagine the
$\TMF$ classes coming from $K$-theory twistings as 1-dimensional elliptic
elements and twisted $K$-theory as $K$-theory with coefficients in
this ``elliptic line bundle''.


\subsubsection{Projective Hilbert Space Bundles}

There is a very simple and well known reformulation of twistings of
$K$-theory as projective Hilbert space bundles and of the corresponding
twisted $K$-theory groups as families of Fredholm operators on these
bundles.  Indeed, the space of unitary operators on Hilbert space is
contractible, so the group of projective unitary operators has the homotopy 
type of $BS^1$.  As such a twisting $\alpha:X \ra K(\Z,3)$ of $K$-theory
determines a projective bundle $\HH(\alpha)$ of Hilbert spaces on $X$.  The 
space of Fredholm operators on a Hilbert space has the homotopy type of 
$\Z \times BU$ and depends only on the projectivization of the Hilbert
space.  Sections of the $\Z \times BU$ bundle associated to the twisting 
$\alpha$ can therefore be thought of as Fredholm operators on the projective
bundle $\HH(\alpha)$.  It remains to develop a general index theory for 
elliptic operators on these projective bundles, but substantial progress has 
been made by Mathai, Melrose, and Singer~\cite{mmr}, who prove an index theorem
in the case that the twisting $\alpha$ is a torsion class in $H^3(X;\Z)$.

\vspace{-7pt}

\subsubsection{$K(\Z,2)$-Equivariant $K$-Theory}

We would like to discuss an algebro-geometric model for twisted
$K$-theory, and the proper formulation is suggested by reinterpreting
twisted $K$-theory as a $K(\Z,2)$-equivariant theory; this formulation will
also hint at connections with the representation theory of
loop groups.  As before, a twisting is a map $\alpha:X \ra BK(\Z,2)$ defining 
a principal $K(\Z,2)$-bundle $P(\alpha)$ on $X$.  The set of sections of the
associated bundle $P(\alpha) \times_{K(\Z,2)} (\Z \times BU)$ is the same
as the set of $K(\Z,2)$-equivariant maps from $P(\alpha)$ to $\Z \times BU$;
that is, the twisted $K$-theory of $X$ is the ``$\! K(\Z,2)$-equivariant'' 
$K$-theory of $P(\alpha)$.  In particular, elements of the twisted $K$-theory
of $X$ are represented by virtual vector bundles on the total space
$P(\alpha)$ of the $K(\Z,2)$-principal bundle associated to the twisting;
these vector bundles $V$ are required to be $K(\Z,2)$-equivariant in the
sense that for a line $L \in K(\Z,2)$, the virtual vector space 
$V_{L\cdot x}$ at the point $L\cdot x \in P(\alpha)$ is equal to
$L \otimes V_x$, for all points $x \in P(\alpha)$.

\vspace*{-7pt}

\subsubsection{Perfect Complexes of $\alpha$-Twisted Sheaves}

Our `space' $X$ will now be a scheme, and a twisting of $K$-theory
is a $\G_m$-gerbe on $X$.  These gerbes are classified by $H^2(X;\G_m)$
and can be thought of as stacks locally isomorphic to the category of
invertible sheaves.  Elements of the twisted $K$-theory of $X$ for a
twisting gerbe $\alpha$ should be virtual sheaves of locally free
$\OO_{\alpha}$-modules on $\alpha$ that are $B\G_m$-equivariant
in an appropriate sense.  More
precisely an element of the twisted $K$-theory of $X$ is a perfect complex
of $\alpha$-twisted sheaves on the gerbe $\alpha$, that is a complex
of $\alpha$-twisted sheaves locally quasiisomorphic to a finite length
complex of free finite rank $\OO_{\alpha}$-modules.  In the topological
situation the analogue of the perfect complex on $\alpha$ is a two term
complex of bundles on $P(\alpha)$, each of countably infinite rank, 
with a differential that is locally an isomorphism off of a finite
rank subbundle.  We would like to emphasize that this notion of 
$\alpha$-twisted $K$-theory elements on the scheme $X$ does not depend on 
the class $\alpha \in H^2(X;\G_m)$ being torsion.

\vspace*{-7pt}

\subsubsection{Central Extensions of Loop Groups}

We now specialize to the case (which indeed will be our primary focus
in this paper) that our space is a connected simply connected compact Lie 
group $G$.  A twisting map $\alpha: G \ra K(\Z,3)$ gives a map from the free
loop space $LG$ to the classifying space $BS^1$ by the composition
$LG \ra LK(\Z,3) \ra \Omega K(\Z,3) \simeq BS^1$, and thereby gives a principal
$S^1$-bundle on $LG$.  The total space $\widetilde{LG}$ of this principal 
bundle can be given a group structure as an $S^1$-central extension of $LG$.
The classifying space $B \widetilde{\Omega G}$ of the based loop central
extension $\widetilde{\Omega G} \subset \widetilde{LG}$ is precisely the
total space $P(\alpha)$ of the principal $K(\Z,2)$ bundle over $G$.  
Moreover, to an irreducible highest-weight representation of
$\widetilde{LG}$ one can associate an equivariant map from $P(\alpha)$
to $\Z \times BU$ and thereby an element of the twisted $K$-theory of
$G$~\cite{mick}.  The precise relation between the representation theory
of loop groups and twisted $K$-theory is described by Freed, Hopkins,
and Teleman~\cite{fht}---they
prove that the group of positive energy unitary representations of
$\widetilde{LG}$ is the twisted $G$-equivariant $K$-theory \mbox{of $G$}.

\subsubsection{B-Fields and D-Branes}

A great deal of the limelight focused on twisted $K$-theory has come from
the widespread realization that certain boundary conditions in string
theory naturally represent elements in the twisted $K$-theory of spacetime.  
In this context the twistings are represented by nontrivial
Neveu-Schwarz B-fields; the elements of twisted $K$-theory are D-branes,
submanifolds of spacetime with a twisted $\Spinc$ structure on
their normal bundles.  More generally, such a submanifold $M$ may be equipped 
with a vector bundle $V$ and the class represented by the pair 
$(M,V)$ is the pushforward of $V$ to the twisted $K$-theory of the ambient 
spacetime $X$.  When the space $X$ is a Lie group, as in this paper, the
twisted $K$-theory can be thought of as a topological model for the space of 
D-branes in a Wess-Zumino-Witten model for conformal field theory.
Frequently the spacetime $X$ is itself a $\Spinc$ manifold;
the D-branes are then twisted $\Spinc$ submanifolds and represent
elements in the twisted $K$-homology of $X$.  In this case, a D-brane $M$
naturally represents a class in a more refined group, the twisted
$\Spinc$ bordism of $X$, and there is a twisted index map
that recovers the twisted $K$-homology class of $M$.  This perspective
guides the discussion of the twisted $\Spinc$ bordism of
Lie groups in the last section of this paper.

\subsection{Results} \label{results}

\vspace{-10pt}

We prove that the twisted $K$-homology ring of a simple Lie group is
an exterior algebra tensor a cyclic group, we give a detailed
description of the orders of these cyclic groups in terms of
the dimensions of irreducible representations of related groups,
and we show that these orders originate, via a twisted index map,
from relations in the twisted $\Spinc$ bordism group. \vspace{8pt}

\begin{theorem} \label{thmtkh}
Let $G$ be a compact, connected, simply connected, simple Lie
group of rank $n$.  The twisted $K$-homology ring of $G$ with
nonzero twisting class $k \in \H^3(G;\Z) \cong \Z$ is an
exterior algebra of rank $n-1$ tensor a cyclic group:
$$K^{\tau(k)}_{\bdot}(G) \cong \Lambda[x_1, \ldots, x_{n-1}] \tensor
\Z/{c(G,k)}. \nn$$
Here $c(G,k)$ is an integer depending on the group and the
twisting.
\end{theorem}
\vspace{7pt}
\nid This fact was first noticed in the case of $SU(n)$ by Hopkins.
The proof is in section~\ref{tatesec} for groups other than
$Spin(n)$, and in section~\ref{torSOn} for $Spin(n)$. 

\begin{theorem} \label{thmcyc}
For the classical groups, the cyclic orders $c(G,k)$, $k>0$, 
of the twisted $K$-homology groups of $G$ are:
\begin{align} \nn
c(SU(n+1),k) &= \gcd{\textstyle \left\{\binom{k+i}{i} - 1 : 
1 \leq i \leq n\right\}} \\
c(Sp(n),k) &= \gcd{\textstyle \left\{\sum\limits_{-k\leq j \leq -1} 
\binom{2j+2(i-1)}{2(i-1)} : 1 \leq i \leq n\right\}} \nn \\
c(Spin(4n-1),k) &= \gcd{\textstyle \left\{\{\binom{k}{i} : 1 \leq i 
\leq 2n-2\} \right. } \nn \\ & \qquad \qquad \qquad \quad {\textstyle \left.
\cup \{2 \binom{k}{2n-1}\} \cup \{2\binom{k}{2i+1}+\binom{k}{2i} :
n \leq i \leq 2n-2\}\right\} } \nn \\
c(Spin(4n+1),k) &= \gcd{\textstyle \left\{\{\binom{k}{i} : 1 \leq i 
\leq 2n-1\} 
\cup \{2 \binom{k}{2i+1}+\binom{k}{2i} : n \leq i \leq 2n-1\}\right\}} \nn \\
c(Spin(4n+2),k) &= \gcd{\textstyle \left\{\{\binom{k}{i} : 1 \leq i
\leq 2n\} \right. } \nn \\ & \qquad\qquad\qquad\quad {\textstyle \left.
\cup \{2\binom{k}{2n+1}\} \cup \{2\binom{k}{2i+1}+\binom{k}{2i}
: n+1 \leq i \leq 2n-1\}\right\} } \nn \\
c(Spin(4n),k) &= \gcd{\textstyle \left\{\{\binom{k}{i} : 1 \leq i
\leq 2n-1\} \cup \{2\binom{k}{2i+1}+\binom{k}{2i} : n \leq i \leq 2n-2
\}\right\}}. \nn 
\end{align}
\end{theorem}
\vspace{-5pt}
\nid (Note that $c(G,-k)=c(G,k)$.  The formulas for $c(Spin(4n-1),k)$
and $c(Spin(4n),k)$ exclude the degenerate case $n=1$.)  
The proofs for $SU(n)$, $Sp(n)$, and
$Spin(n)$ occur respectively in sections~\ref{genSUnG2},~\ref{genSpn},
and~\ref{torSOn}.  A general method for computation, applicable
to the exceptional groups, is discussed in section~\ref{cycgen},
and the cyclic order for $G_2$ is given in section~\ref{genSUnG2}.

The referee has drawn our attention to an intriguing conjecture by Volker Braun
concerning these cyclic orders~\cite{braun}, namely that $c(G,k)$ is the greatest
common divisor of the dimensions of the representations generating the
Verlinde ideal.  (The Verlinde ideal is the ideal $I$ in the representation ring $R[G]$
such that the quotient $R[G]/I$ is the Verlinde algebra; by Freed, Hopkins, and 
Teleman~\cite{fht} the Verlinde algebra is the twisted $G$-equivariant $K$-theory
of $G$.)  Using a Kunneth spectral sequence in equivariant $K$-theory due to Hopkins, 
Braun is able to show that
the conjecture holds for any group for which the Verlinde algebra is a complete
intersection---for such groups $G$ Braun's argument also shows that the twisted $K$-theory of $G$ is a
group isomorphic to the exterior algebra described in Theorem~\ref{thmtkh} above; Braun does not address the question of which groups satisfy the complete intersection
condition.  Using a twisted Bousfield-Kan spectral sequence one can see that
for the groups $SU(n)$ and $Sp(n)$ the Verlinde algebra is a complete
intersection~\cite{cldverl}.  Computational evidence strongly suggests that Braun's conjecture
is true for all groups, but at the moment it remains open for the spin and for the exceptional
groups.

\vspace{8pt}

\begin{propos}
Let $G$ be as in Theorem~\ref{thmtkh}.  Suppose $M_i$ is a collection of
$\Spinc$ manifolds over $\Omega G$ whose fundamental classes
generate $K.\Omega G$ as an algebra.  Then there are twisted
$\Spinc$ structures on the bordisms $W_i = M_i \times I$ such 
that the cyclic order of the twisted $K$-homology of $G$ is
$\gcd(\ind(\partial W_1), \ldots, \ind(\partial W_n))$,
where $\ind: MSpin^{c}_{\bdot} \ast \ra K. \ast$ is the
index map from $\Spinc$ bordism to $K$-homology.
\end{propos}
\vspace{2pt}
\nid The proof of this proposition is the focus of section~\ref{nullbord}.

\subsection{Techniques and Overview}

\vspace{-5pt}

The primary tool for calculating twisted $K$-homology rings is the
twisted Rothenberg-Steenrod spectral sequence; this is the original
method used by Hopkins in the case $G=SU(n)$.  The spectral
sequence is:
$$E^2 = \Tor^{K.\Omega G}(\Z,\Z_{\tau(k)}) \Rightarrow
K^{\tau(k)}_{\bdot}(G),$$
where $\Z_{\tau(k)}$ is the integers with a twisted 
$K.\Omega G$-module structure depending on $k$.
In section~\ref{tht} we present various generalities about twisted
homology theories; then in section~\ref{trsss} we use a method
of Segal~\cite{segalcsss} to construct this Rothenberg-Steenrod
spectral sequence in twisted $K$-homology.

As the $K$-homology rings of loop spaces of simple Lie groups are known,
our primary task is computing the $\Tor$ groups over these rings.
Remarkably, for $G \neq Spin(n)$ this can be done without identifying
the twisted $K.\Omega G$-module structure on $\Z$.  These $\Tor$
groups are calculated in section~\ref{tatesec} by an iterated series
of filtration spectral sequences applied to a judiciously chosen
Tate resolution.  The spectral sequences are seen to collapse and to be
extension-free, completing the proof of Theorem~\ref{thmtkh} 
for $G \neq Spin(n)$.

The $\Tor$ computation for $\Spin(n)$ requires a detailed knowledge
of the twisted module structure on $\Z$; this module structure is
also precisely what is needed to identify the cyclic orders of
the twisted $K$-homology groups.  The best way to identify this
module structure is via generating varieties for the loop space
of the group, and this is the subject of section~\ref{gvco}.
Sections~\ref{genSUnG2},~\ref{genSpn}, and~\ref{cycgen} describe
generating varieties for various groups, compute the cyclic orders
in the corresponding cases, and discuss a general method for
determining the cyclic order.  Section~\ref{torSOn} describes the 
twisted module structure for $Spin(n)$ and presents the belated $\Tor$
calculation for this group.

The computation in section~\ref{gvco} of the cyclic order in terms of 
the dimensions of irreducible representations does not
give much geometric insight into these torsion groups.  We give,
in section~\ref{bord}, an interpretation of these orders in terms 
of relations in the twisted $\Spinc$ bordism group of $G$.
The main tool, presented in section~\ref{cocycle}, is a
cocycle model for twisted $\Spinc$ bordism.  This model
allows explicit descriptions of nullbordisms of particular
$\Spinc$ manifolds over $G$ corresponding to relations in the
twisted $K$-homology of $G$---see section~\ref{nullbord}.  We conclude in 
section~\ref{exterior} by discussing potential representatives in
$MSpin^{c,\tau}_{\bdot}(G)$ for the exterior generators of 
$K^{\tau}_{\bdot}(G)$.

\vspace{-5pt}

\section{Twisted $K$-Theory and the Rothenberg-Steenrod Spectral Sequence}

\vspace{-5pt}

\subsection{Twisted Homology Theories} \label{tht}

\vspace{-5pt}

We review the definitions and basic properties of twisted homology and
cohomology theories.  There are by now various models for these theories,
but the following perspective owes as much to Goodwillie as to folklore.

For a spectrum $F$, the cohomology of a space $X$ with coefficients in
$F$ can be defined as
$$F^n(X) := \colim \Gamma_h(X, X \times \Omega^i F_{i+n}); \nn $$
here $\Gamma_h(X,E)$ refers to homotopy classes of sections of the
(here trivial) bundle $E$ on $X$.  The maps in the colimit are induced by 
applying the usual structure maps $\Omega^i F_{i+n} \ra \Omega^{i+1} 
\Sigma F_{i+n} \ra \Omega^{i+1} F_{i+1+n}$ fibrewise to the bundle
$X \times \Omega^{i} F_{i+n} \ra X$.  Now let $E$ be a bundle of based
spectra over $X$, with fibre spectrum $F$; this means in particular
that for each $i$ we have a fibration $E_i \ra X$, a section
$X \ra E_i$, and a fibrewise structure map $\Sigma_X E_i \ra E_{i+1}$.
(Note that $\Sigma_X$ denotes fibrewise suspension and $\Omega_X$ will
denote the fibrewise loops.)  The cohomology of $X$ with coefficients 
in $E$ is defined to be
$$E^n(X) := \colim \Gamma_h (X, \Omega_X^i E_{i+n}) \nn $$
where the colimit maps are, as expected, induced by
$\Omega_X^i E_{i+n} \ra \Omega_X^{i+1} \Sigma_X E_{i+n} \ra
\Omega_X^{i+1} E_{i+1+n}$.

The parallel in homology is similar.  The homology of $X$ with
coefficients in $F$ is
$$F_n(X) := \colim [S^{i+n}, (X \times F_i) / X], \nn $$
with maps induced by $\Sigma ((X \times F_i)/X) =
(X \times \Sigma F_i) / X \ra (X \times F_{i+1}) / X.$
As above, when $E$ is a bundle of based spectra, we have a `base point'
section $X \ra E_i$ for all $i$.  The homology of $X$ with
coefficients in $E$ is
$$E_n(X) := \colim [S^{i+n}, E_i / X]; \nn $$
the colimit maps are induced by $\Sigma (E_i / X) = (\Sigma_X E_i) / X
\ra E_{i+1} / X.$

For completeness we also mention the reduced analogs of homology and
cohomology with coefficients in a bundle of spectra.  The reduced
cohomology with coefficients in a trivial $F$ bundle can be given
as
$$\tilde{F}^n(X) := \colim \Gamma_h^b(X, X \times \Omega^i F_{i+n}), \nn $$
that is as the colimit of homotopy classes of sections taking
the base point of $X$ to the basepoint of $\Omega^i F_{i+n}$.  The
reduced cohomology with coefficients in $E$ is then
$$\tilde{E}^n(X) := \colim \Gamma_h^b (X, \Omega_X^i E_{i+n}); \nn $$
the maps are induced as before.  Similarly, the reduced homology
with coefficients in a trivial bundle is
$$\tilde{F}_n(X) := \colim [S^{i+n}, (X \times F_i) / (X \vee F_i)]. \nn $$
The twisted reduced homology is finally
$$\tilde{E}_n(X) := \colim [S^{i+n}, E_i / (X \vee F_i)]; \nn $$
the maps are induced by $\Sigma (E_i / (X \vee F_i))
= (\Sigma_X E_i) / (X \vee \Sigma F_i) \ra E_{i+1} / (X \vee F_{i+1}).$
Of course, these reduced groups are special cases of the relative groups:
$$E^n(X,A) := \colim \Gamma_h(X,A; \Omega_X^i E_{i+n}, s(A)), \nn $$
where $s$ is the distinguished base point section; similarly,
$$E_n(X,A) := \colim [S^{i+n}, (E_i / X) / ((E_i|_A) / A)]. \nn $$

The most important fact about twisted homology theories is that they
are honest homology theories in an appropriate category.  Indeed, consider
the category of pairs $(X,A)$ of spaces, where $A$ is a closed subspace
of $X$ and $X$ is equipped with a bundle $E$ of based spectra with fibre
spectrum $F$.  From the above description of the homology $E_n(X,A)$,
it is immediate that twisted homology on this category of pairs is a
homology theory in the classical sense.

\begin{center}$\sim\sim\sim\sim\sim$\end{center} \vspace{6pt}

In this paper we will only be concerned with bundles of spectra associated
to principal $K(\Z,2)$ bundles over our space $X$.  As usual, we fix
a model for $K(\Z,3)$ and select a particular universal $K(\Z,2)$ bundle
on it.  A map $\alpha: X \ra K(\Z,3)$ gives a principal
$K(\Z,2)$ bundle $P(\alpha)$ on $X$, classified up to isomorphism by the 
homotopy class of the map.  For any basepoint-preserving action of
$K(\Z,2)$ on a spectrum $F$, we can form the associated $F$ bundle to
$P(\alpha)$.  The resulting bundle $P(\alpha) \times_{K(\Z,2)} F$ is a
bundle of based spectra on $X$, as above.  Note that on the level of
spaces, the action of $K(\Z,2)$ on $F$ is given by maps
$K(\Z,2)_+ \sm F_i = (K(\Z,2) \times F_i)/(K(\Z,2) \times *) \ra F_i$,
and we often denote the spectrum action simply by a map $K(\Z,2)_+ \sm F
\ra F$.

Our primary examples are twisted $\Spinc$-bordism and twisted $K$-theory.
The $K(\Z,2)$ bundle
$$ K(\Z,2) = BU(1) \ra B \Spinc \ra BSO \nn $$
is principal, with classifying map $BSO \xra{\beta w_2} BBU(1) = K(\Z,3)$
classifying the integral Bockstein of the second Stiefel-Whitney class.
In particular we have an action $K(\Z,2) \times B\Spinc \ra B\Spinc$;
on Thom spaces this action is $K(\Z,2)_+ \sm M\Spinc \ra M\Spinc$,
that is, a based action of $K(\Z,2)$ on the $\Spinc$ Thom spectrum.
The $\alpha$-twisted $\Spinc$-bordism groups are then, of course, the
stable homotopy groups $\pi_i((P(\alpha) \times_{K(\Z,2)} M\Spinc)/X)$.

The $K$-theory spectrum $K$ is a module over $\Spinc$-bordism by the
usual index map $M\Spinc \xra{\ind} K$.  Taking the above based action
$K(\Z,2)_+ \sm M\Spinc \xra{\phi} M\Spinc$ and smashing over $M\Spinc$
with $K$, we have a compatible based action on $K$-theory:
\begin{equation} \nn
\xymatrix@C+5pt{
K(\Z,2)_+ \sm M\Spinc \ar@{->}[r]_-{\phi} \ar@{->}[d]_{\id \sm \ind} 
& M\Spinc \ar@{->}[d]^{\ind}\\
K(\Z,2)_+ \sm K \ar@{->}[r]_-{\phi \sm_{M\Spinc} (\id)} & K
}
\end{equation}
The corresponding map on associated principal bundles
$P(\alpha) \times_{K(\Z,2)} M\Spinc \ra P(\alpha) \times_{K(\Z,2)} K$
induces a map from twisted $\Spinc$-bordism to twisted $K$-theory
which we call the twisted index map.  This map will be important
in section~\ref{bord}.

Twisted $K$-theory can be defined more directly by choosing an
explicit model for $\Z \times BU$ (typically the space of Fredholm
operators on a fixed Hilbert space $\HH$) that admits an explicit
action by some model for $BU(1)$ (typically the space of projective
unitary operators on $\HH$); see, for example, Atiyah~\cite{atiyahkpp}.
(The referee has pointed out that the recent exposition by Atiyah and 
Segal~\cite{astkt} is another source for operator-theoretic definitions
of twisted $K$-theory.)
Whatever the formal definition, the geometric action being modeled is 
the following: a complex line $L$ (representing a point in $BU(1)$) acts 
on a virtual-dimension-zero (or stable) vector space $V$ (representing a 
point in $BU$) by tensor product, that is, $V \mapsto L \tensor V$.

It is worth noting, though, that this heuristic action of tensoring a vector
bundle with a line can be misleading
if we pay insufficient attention to the virtual dimension zero condition.
It is tempting to think of elements of $\alpha$-twisted $K$-cohomology as
sections of an $\alpha$-twisted gerbe of rank $n$, for some sufficiently
large $n$; (such a section is locally a rank-$n$ vector bundle, twisted
globally by $\alpha$).  However, in this paper we are dealing with 
non-torsion twistings, and therefore no nontrivial element of twisted 
$K$-cohomology is representable by a section of any finite rank gerbe. 
We are inescapably in either a virtual-dimension-zero or an 
infinite-dimensional situation---which would seem to be a matter of 
personal penchant.

\subsection{The Twisted Rothenberg-Steenrod Spectral Sequence} \label{trsss}

The ``twisted'' Rothenberg-Steenrod spectral sequence computing the
twisted $K$-homology of a space is in fact the ordinary Rothenberg-Steenrod
(a.k.a. homology Eilenberg-Moore) spectral sequence in an 
appropriate category, and as such requires little comment.  We briefly
recall the spectral sequence in generality, then describe its application 
to the geometric bar complex on the loop space of a simple Lie group.

We work in the category $\mathcal{K}$ of pairs $(X;E)$, where $X$ is a space 
and $E$ is a bundle of based spectra on $X$ with fibre the $K$-theory spectrum;
the morphisms are those bundle maps that are homotopy equivalences on each 
fibre. Similarly, we have a category of triples $(X,A;E)$ where $A$ is a closed
subspace of $X$ and $E$ is again a bundle on $X$.  As mentioned in the
last section, the functors $$(X,A;E) \mapsto E_n(X,A) = \colim [S^{i+n}, 
(E_i / X) / ((E_i|_A) / A)]\nn$$ form a homology theory in the classical sense.
In particular, for any simplicial object $S.$ in $\mathcal{K}$, there is a 
spectral sequence a la Segal~\cite{segalcsss} with $E^2$ term $H_p(E_q(S.))$
converging to the homology of the realization $E_{p+q}(|S.|)$.

Let $G$ be a simple, simply connected Lie group and 
$k \in H^2(\Omega G;\Z) = \Z$ an integer
describing a line bundle $L^{-k}$ on the loop space $\Omega G$.  On the one 
hand, there is the trivial projection map in $\mathcal{K}$ from $(\Omega G;
\Omega G \times K)$ to $(\ast;K)$.  On the other hand, there is a twisted map
$\tau(k):(\Omega G;\Omega G \times K) \ra (\ast;K)$ given by
$\Omega G \times K \overset{k \times \id}{\longrightarrow} K(\Z,2) \times K 
\ra K$, where the last map is the $K(\Z,2)$ action on the spectrum $K$
described in section~\ref{tht}.  The geometric bar construction 
$B_{\tau}\Omega G = B.(\ast, \Omega G, \ast_{\tau})$ is a simplicial object 
in $\mathcal{K}$.
To describe the corresponding spectral sequence we need only compute the
effect of $\tau(k)$ in homology and identify the realization 
$|B_{\tau}\Omega G|$.  

Given a class $\phi$ in the $K$-homology of $\Omega G$ the
image of $\phi$ under $\tau(k)$ is evidently equal to the evaluation
$\langle\tau(k)^*(1),\phi\rangle$, where $\langle -,- \rangle$ denotes the 
Kronecker pairing.  The pullback $\tau(k)^*(1)$ is $L^{k}$, and the resulting 
map $K.\Omega G \overset{\langle L^k,-\rangle}{\longrightarrow} K.\ast$ 
defines a module structure on $K.\ast$ which
we denote $(K.\ast)_{\tau}$.  The $E^2$ term of our spectral sequence is
therefore $\Tor^{K.\Omega G}(K.\ast,(K.\ast)_{\tau})$.  

As a space the
realization of $B_{\tau}\Omega G$ is evidently $B\Omega G \simeq G$; we
identify the $K$-bundle.  The $K$-bundle on the realization is defined
by a 1-cocycle $\tau(k)$ with values in $K(\Z,2)$ and as such is classified 
by the image of $\tau(k)$ in $H^3(B\Omega G;\Z)$.  We have
$H^3(B\Omega G) \cong H^3(\Sigma \Omega G)$ and it is enough to identify 
the restriction of $\tau(k)$ to the 1-skeleton $\Sigma \Omega G$ of 
$B\Omega G$.  It is, however, immediate that this cocycle on the 1-skeleton of
the geometric bar construction $B_{\tau}\Omega G$ has homology invariant
$k \in H^3(\Sigma \Omega G)$.  In summary: \vspace{8pt}
\begin{propos}
There is a spectral sequence of algebras with $E^2$ term
$$E^2_{pq} = \Tor^{K.\Omega G}_{p,q}(K.\ast,(K.\ast)_{\tau}) \nn$$
converging as an algebra to the twisted $K$-homology
$K_{p+q}^{\tau} (G)$.
\end{propos}
\vspace{5pt}

The twisted $K$-homology of $G$ is by definition the homotopy of the spectrum
$E/G$, where $E$ is the bundle of spectra over $G$ determined by the
twisting class $\tau$.  There is a pairing of spectra $E/G \sm E/G \ra E/G$
induced by the multiplication on $G$, and this pairing gives the algebra
structure on $K_{\bdot}^{\tau}(G)$.  Note that the existence of this pairing
depends essentially on the fact that the twisting class $\tau$ is primitive as
an element of $H^3(G)$; because $G$ is simple and simply connected, all 
such elements are indeed primitive.

The multiplicative structure of the above spectral sequence
can be seen as follows.  The filtration of $B \Omega G \simeq G$ by the standard 
skeleta $B_i \Omega G$ induces
a filtration of $E/G$ by a tower $T$ of spectra $T_i = (E |_{B_i \Omega G}) / B_i \Omega G$.
The spectral sequence in question is the homotopy spectral sequence associated to this
tower.  The multiplication on $G$ corresponds to a filtration-preserving multiplication on
$B \Omega G$.
The pairing on $E/G$ therefore induces a pairing of towers $T \sm T \ra T$ and this pairing
of towers descends to the algebra structure on the homotopy spectral sequence---see for
example the careful exposition of homotopy spectral sequence pairings by 
Dugger~\cite{dugger}.  That the resulting algebra structure on the $E^2$ term of the
spectral sequence agrees with the usual pairing on $\Tor$ follows by comparing the
algebraic and geometric bar constructions for $K.\Omega G$ and $\Omega G$ respectively, as
in for instance~\cite{rs,may,rao}.

\section{Tate Resolutions and $\Tor^{K.\Omega G}(\Z,\Z_{\tau})$
for $G \neq Spin(n)$} \label{tatesec}

For each group $G$, we describe the $K$-homology of the loop space of 
$G$, give an appropriate Tate resolution of $K.\ast=\Z$ over $K.\Omega G$,
and compute the torsion group using a series of filtration
spectral sequences.

We recall Tate's main result on algebra resolutions over a commutative
Noetherian ring $R$.  An ideal $I \subset R$ is said to be generated
by the regular sequence $a_1, \ldots, a_r \in R$ if $I=(a_1, \ldots,
a_r)$ and $a_i$ is not a zero-divisor in $R/(a_1, \ldots, a_{i-1})$
for all $i$. \vspace{8pt}
\begin{theorem}[Tate~\cite{tate}]
Let $A \subset B$ be ideals of $R$ generated respectively by
the regular sequences $(s_1, \ldots, s_m)$ and $(t_1, \ldots, t_n)$.
For any choice of constants $c_{ji} \in R$ such that
$s_j = \sum_{i=1}^n c_{ji} t_i$, the differential graded algebra
$$D=(R/A \langle T_1, \ldots, T_n \rangle \{S_1, \ldots, S_m\}; d(T_i) = [t_i],
d(S_j) = \sum_{i=1}^n [c_{ji}] T_i) \nn $$
is a resolution of $R/B$ as an $R/A$-module.  Here the $T_i$ are
strictly skew commutative generators of degree 1, and the $S_j$
are divided power algebra generators of degree 2.
\end{theorem}
\vspace{10pt}
\nid In particular, $\Tor_{R/A}(R/B,Q)$ will be given as the
homology $H(D \tensor_{R/A} Q)$.  In our applications, 
$R$ will be a polynomial ring $\Z[x_1, \ldots, x_n]$, the
ideal $A$ will depend on the group, the ideal $B$ will be
$(x_1, \ldots, x_n)$, and $Q$ will be an $R/A$-module
$\Z_{\tau}$ on which $x_i$ acts by an integer $c_i$ depending on
the group and the twisting class.

\subsection{$\Tor$ for $SU(n+1)$ and $Sp(n)$} \label{torSUnSpn}

Elementary calculation shows that the integral cohomology rings
of $SU(n+1)$ and $Sp(n)$ are exterior algebras on $n$ generators.
Application of the spectral sequence $\Ext^{H^{\bdot}(G;k)}(k,k) 
\Rightarrow H.(\Omega G; k)$, $k$ a field, then implies that
the integral Pontryagin rings $H.(\Omega SU(n+1))$ and 
$H.(\Omega Sp(n))$ are both polynomial on $n$ generators, all
in even degree.  In each case the Atiyah-Hirzebruch spectral 
sequence for $K$-theory then collapses, and the $K$-theory Pontryagin 
ring is again polynomial.

The Tate resolution in this case is especially simple, as
the ideal $A$ is trivial.  Let $G$ denote either $SU(n+1)$ or
$Sp(n)$ and $k \in \Z \cong H^3(G;\Z)$ the twisting class.
Choose reduced generators $x_i$ of $K.\Omega G$, so that $K.\Omega G
\cong \Z[x_1, \ldots, x_n]$.  (Note that, unless otherwise 
noted, we treat $K$-theory as $\Z/2$-graded.)  
The $K.\Omega G$ module 
structure on $\Z_{\tau}$ is given, as in 
section~\ref{trsss}, by the map $K.\Omega G \ra K.*$ sending
a class $x$ to $\langle L^k, x \rangle$,
where $L$ is a generating line bundle.  We defer the
explicit evaluation of these maps to section~\ref{gvco}.
For now, we denote by $c_i$ the image of $x_i$ in $\Z_{\tau}$;
of course this constant depends on both the group and the
twisting, but we tend to omit both dependencies from the notation.
By Tate's theorem,
\begin{align} 
\Tor^{K.\Omega G}(\Z,\Z_{\tau}) &= 
H(\Z[x_1, \ldots, x_n]\langle T_1, \ldots, T_n \rangle 
\tensor_{\Z[x_1, \ldots,
x_n]} \Z_{\tau}; d) \nn \\
&= H(\Z\langle T_1, \ldots, T_n \rangle; dT_i = c_i). \nn
\end{align}
\nid To evaluate this homology group we employ the following general
procedure.  Suppose we know the homology of the subalgebra generated
by $T_1, \ldots, T_i$.  We filter the subalgebra generated
by $T_1, \ldots, T_{i+1}$ by powers of $T_{i+1}$ and look at the
associated spectral sequence.  The only differential is $d^1$, which
is given by multiplication by $c_{i+1}$, and by induction we can thus
compute the homology of the original algebra.  

We assume for now that $c_1$ is not zero; this is indeed the case (see 
sections~\ref{genSUnG2} and~\ref{genSpn}).  The homology of $(\Z
\langle T_1 \rangle,d)$
is $\Z/c_1$.  The, quite degenerate, spectral sequence of the filtration
of $(\Z \langle T_1, T_2 \rangle,d)$ by $T_2$ is therefore
$$\Z/c_1 \xla{c_2} \Z/c_1. \nn $$
The homology is $\Z/g_{12}\langle y_2 \rangle$, where $g_{12} = 
\gcd\{c_1,c_2\}$ and
$y_2$ is an exterior class.  More generally we will denote by
$g_{1..i}$ the greatest common divisor $\gcd\{c_1,c_2, \ldots,
c_i\}$.  The induction step is, as expected, the homology of
$$\Z/g_{1..i}\langle y_2, \ldots, y_i \rangle \xla{c_{i+1}} 
\Z/g_{1..i} \langle y_2, \ldots, y_i \rangle, \nn $$
and the $\Tor$ groups are given by
\begin{align}
\Tor^{K.\Omega SU(n+1)}(\Z,\Z_{\tau}) &= \Z/(g_{1..n}(SU(n+1),k))
\langle y_2, \ldots, y_{n-1} \rangle \nn \\
\Tor^{K.\Omega Sp(n)}(\Z,\Z_{\tau}) &= \Z/(g_{1..n}(Sp(n),k))
\langle y_2, \ldots, y_{n-1} \rangle. \nn
\end{align}
We belabor this calculation only because, when we come to
more complicated examples, especially $Spin(n)$, 
it will help to have a clear model.

\vspace{-4pt}
\subsection{$\Tor$ for the Exceptional Groups} \label{torG2F4E6}

\vspace{-6pt}
The exceptional Lie groups are nature's best attempts to make
a finite dimensional Lie group out of $K(\Z,3)$.  In particular
they are homotopy equivalent to $K(\Z,3)$ through a range of
dimensions, and so their loop spaces are homotopy equivalent
to $K(\Z,2)$ through a similar range.  The $K$-homology of
$K(\Z,2)$ is the subalgebra of $\Q[a]$ generated by
$\{a, \binom{a}{2}, \binom{a}{3}, \ldots\}$; see~\cite{adams}.  
Extensive computations
by Duckworth~\cite{duck} show that for $G$ exceptional,
the $K$-homology $K.\Omega G$ differs from a polynomial ring only in the 
aforementioned low-dimensional flirtation with $K(\Z,2)$.  For example,
Duckworth proves that $K.\Omega E_8$ is a polynomial ring on seven
generators tensor the subalgebra of $\Q[a]$ generated by the elements
$\{a, \binom{a}{2}, \binom{a}{3}, \binom{a}{4}, \binom{a}{5}\}$.  In
order to use Tate resolutions, we must give explicit algebra presentations
of these $K$-homology rings: \vspace{8pt}
\begin{propos} \label{proppres}
The $K$-homology rings of the loop spaces of the 
exceptional Lie groups are given by
\vspace{-6pt}
\begin{align}
K.\Omega G_2 &= \frac{\Z[a,b,x_3]}{(a(a-1)-2b)} \nn \\
K.\Omega F_4 &= \frac{\Z[a,b,c,x_4,x_5,x_6]}{(a(a-1)-2b, b(a-2)-3c)} \nn \\
K.\Omega E_6 &= \frac{\Z[a,b,c,x_4,x_5,x_6,x_7,x_8]}
{(a(a-1)-2b, b(a-2)-3c)} \nn \\
K.\Omega E_7 &= \frac{\Z[a,b,c,d,x_5,x_6,x_7,x_8,x_9,x_{10}]}
{(a(a-1)-2b, b(a-2)-3c, b(b+1) - a(b+c) - 2d)} \nn \\
K.\Omega E_8 &= \frac{\Z[a,b,c,d,e,x_6,x_7,x_8,x_9,x_{10},x_{11},x_{12}]}
{(a(a-1)-2b, b(a-2)-3c, b(b+1) - a(b+c) - 2d, d(a-4)-5e)} \nn
\end{align}
\end{propos}
\nid Note that the unsightly third relation in the rings for $E_7$ and $E_8$
is essential and cannot be replaced by the more sensible relation
$c(a-3)-4d$.
We remark that, because the `lettered' generators in these $K$-homology rings
come from corresponding generators in $K.(K(\Z,2))$, the twisted
pushforwards of these elements are easily computed.  In particular,
the twisted pushforward of $a$, denoted again by $c_1$, is just $k$,
the twisted pushforward of $b$ is $c_2=\binom{k}{2}$, of $c$ is
$c_3=\binom{k}{3}$, and so on, with each generator mapping to its
respective binomial coefficient.

As always, our starting point is the Tate resolution:
$$\Tor^{K.\Omega G_2} (\Z,\Z_{\tau}) =
H(\Z\langle T_1, T_2, T_3\rangle\{S_1\}; 
dT_i=c_i, dS_1=(c_1 - 1)T_1 - 2T_2). \nn $$
Consider the subalgebra generated by $T_1$, $T_2$, and $S_1$.
If $k$ is even, we can rewrite this DGA as
$$ (\Z\langle T_1',T_2'\rangle\{S_1\}; dT_1'=0, dT_2'=\frac{k}{2}, 
dS_1=T_1'), \nn $$
where $T_1'=(k-1)T_1 - 2T_2$ and $T_2'=\frac{k}{2}T_1 - T_2$.
The Kunneth theorem immediately shows that the homology of this
DGA is $\Z/(\frac{k}{2})$.  If $k$ is odd, we instead change the basis to
$T_1'=\frac{k-1}{2} T_1 - T_2$ and $T_2'=kT_1 - 2T_2$.  The algebra then
takes the form
$$ (\Z\langle T_1',T_2' \rangle\{S_1\}; dT_1'=0, dT_2'=k, dS_1=2 T_1'), \nn $$
and by the Kunneth theorem its homology is $\Z/k$.  In other words,
the homology of the subalgebra in question is, in any case,
$\Z/g_{12}$, where as before $g_{12}=\gcd\{c_1,c_2\}$.  Filtering
as in section~\ref{torSUnSpn} we see that the full $\Tor$ group
is $\Z/g_{123}\langle y_3 \rangle$.

The Tate resolution for $F_4$ gives
\begin{multline} \nn
\Tor^{K.\Omega F_4} (\Z,\Z_{\tau}) =
H(\Z\langle T_1,T_2,T_3,T_4,T_5,T_6 \rangle\{S_1,S_2\}; \\
dT_i=c_i, dS_1=(c_1 - 1)T_1 - 2T_2, dS_2=(c_1 - 2)T_2 - 3T_3). \nn
\end{multline}
\nid We focus on the subalgebra generated by $\{T_1,T_2,T_3,S_1,S_2\}$.
The method used for $G_2$, of changing basis to split the algebra 
into simpler pieces, works here as well; the basis change now depends
on $k$ modulo 6.  We spell out only the case $k=1 \,({\rm mod}\, 6)$.  As
basis change for the $T_i$'s we take
\begin{displaymath} \nn
\left( \begin{array}{ccc}
\frac{k-1}{2} & -1 & 0 \\
-\frac{k-1}{6} & \frac{k-1}{3} & -1 \\
\frac{3k-1}{2} & -1-k & 3
\end{array} \right).
\end{displaymath}
The algebra then has the form $$(\Z\langle T_1',T_2',T_3'\rangle\{S_1,S_2\};
dT_1=dT_2=0, dT_3=k, dS_1=2T_1', dS_2=3T_2'+T_1'). \nn$$  

\pagebreak

\nid The spectral sequence associated to the filtration of the
$\{T_1',T_2',S_1,S_2\}$ subalgebra by powers of $S_2$ is
\begin{equation} \nn
\xymatrix@C=20pt@R=3pt{
*+\txt{$\Z/2$ \\ $\sss (S_1^{(2)}T_1'T_2')$} & 
*+\txt{$\Z/2$ \\ $\sss (S_2S_1^{(2)}T_1')$} \ar[l] & 
*+\txt{$\Z/2$ \\ $\sss (S_2^{(2)}S_1T_1'T_2')$} & \cdots \ar[l] \\
*+\txt{$\Z/2$ \\ $\sss (S_1^{(2)}T_1')$} & 
*+\txt{$\Z/2$ \\ $\sss (S_2S_1T_1'T_2')$} & 
*+\txt{$\Z/2$ \\ $\sss (S_2^{(2)}S_1T_1')$} \ar[l] & \cdots \\
*+\txt{$\Z/2$ \\ $\sss (S_1T_1'T_2')$} & 
*+\txt{$\Z/2$ \\ $\sss (S_2S_1T_1')$} \ar[l] & 
*+\txt{$\Z/2$ \\ $\sss (S_2^{(2)}T_1'T_2')$} & \cdots \ar[l] \\
*+\txt{$\Z/2$ \\ $\sss (S_1T_1')$} & 
*+\txt{$\Z/2$ \\ $\sss (S_2T_1'T_2')$} & 
*+\txt{$\Z/2 \oplus \Z$ \\ $\sss (S_2^{(2)}T_1',S_2^{(2)}T_2')$} 
\ar[l]_-{(1,1)} & \cdots \ar[l]\\
*+\txt{$\Z/2$ \\ $\sss (T_1'T_2')$} & 
*+\txt{$\Z/2 \oplus \Z$ \\ $\sss (S_2T_1',S_2T_2')$} \ar[l]_-{(1,1)} & 
*+\txt{$\Z$ \\ $\sss (S_2^{(2)})$} \ar[l]_-{(1,3)} & \\
*+\txt{$\Z/2 \oplus \Z$ \\ $\sss (T_1',T_2')$} & 
*+\txt{$\Z$ \\ $\sss (S_2)$} \ar[l]_-{(1,3)} & & \\
*+\txt{$\Z$ \\ $\sss (1)$} 
\save[]+<-1cm,-.6cm> \ar@{-}[uuuuuu]+<-1cm,.5cm>
\ar@{-}[rrr]+<.3cm,-.6cm>
\restore
& & & 
}
\end{equation}
There are, of course, no differentials beyond $d^1$ and the homology of 
the $\{T_1',T_2',S_1,S_2\}$ subalgebra is therefore $\Z/6$ in odd degree, 
$0$ in positive even degree, and $\Z$ in degree zero; consequently the 
homology of the $\{T_1',T_2',T_3',S_1,S_2\}$ subalgebra is $\Z/k$
concentrated in degree zero. In general, ie for $k$ not necessarily
congruent to 1 modulo 6, this $\Z/k$ is replaced by $\Z/g_{123}$ and the full 
$\Tor$ group for $F_4$ is $\Z/g_{1..6}\langle y_4,y_5,y_6\rangle$.  
The computation for $E_6$ is
identical, but for two additional exterior generators in the final 
$\Tor$ group.

This basis change approach quickly becomes impractical: for $E_8$ the
congruence of $k$ modulo 60 determines the structure of the basis change
and of the subsequent homology computation.  If we are willing to give
up our ability to write down explicit generators for the $\Tor$ groups,
we can do the computation without such a case by case analysis.
We briefly reconsider the groups $G_2$ and $F_4$.  For $G_2$ the main step
was computing the homology of the DGA $$D=(\Z\langle T_1, T_2\rangle\{S_1\}; 
dT_i=c_i, 
dS_1=(c_1 - 1)T_1 - 2T_2);\nn$$ recall that $c_1=k$ and $c_2=\binom{k}{2}$.  
The homology of the $\{T_1,T_2\}$ subalgebra is $\Z/g_{12}\langle y_2\rangle$, 
where the generator $y_2$ can be taken to be $-(c_1/g_{12})T_2$ 
modulo terms involving $T_1$.  (We will refer to terms with lower 
indices, sensibly enough, as `lower terms' and so say, for example, 
that ``$y_2$ is $-(c_1/g_{12})T_2$ modulo lower terms'').  Thus, 
when we filter $D$ by powers of $S_1$, the homology of $D$ becomes the
homology of
$$(\Z/g_{12}\langle y_2\rangle\{S_1\}; 
dS_1=(2g_{12}/c_1)y_2).\nn$$
Note that $2g_{12}/c_1$ is an integer, so this expression makes sense.
We observe that $2g_{12}/c_1$ is actually a unit in $\Z/g_{12}$;
indeed $g_{12} = c_1 /\!\gcd(2,c_1)$ so
\vspace{-2pt}
$$\gcd(2g_{12}/c_1,g_{12}) = \gcd(2/\!\gcd(2,c_1),c_1/\!\gcd(2,c_1)) = 1. \nn $$
The homology of $D$ is therefore simply $\Z/g_{12}$, as previously noted,
and thus the full $\Tor$ group is again an exterior algebra tensor a cyclic
group.

The case of $F_4$ (and therefore of $E_6$) is again similar.  The main
step is the computation of the homology of the DGA 
\vspace{-2pt}
$$D=(\Z\langle T_1,T_2,T_3\rangle\{S_1,S_2\};
dT_i=c_i, dS_1=(c_1 - 1)T_1 - 2T_2, dS_2=(c_1 - 2)T_2 - 3T_3).\nn$$  (Here
again $c_i=\binom{k}{i}$.)  Using the $G_2$ result we see that
the homology of the $\{T_1,T_2,T_3,S_1\}$ subalgebra is
$\Z/g_{123}\langle y_3 \rangle$ where $y_3$ is $-(g_{12}/g_{123}) T_3$ 
modulo lower terms.
As above the homology of $D$ is thereby reduced to the homology of
\vspace{-2pt}
$$(\Z/g_{123}\langle y_3 \rangle\{S_2\}; dS_2=(3g_{123}/g_{12})y_3).\nn$$
Again, this differential is an isomorphism, ie $3g_{123}/g_{12}$ is a
unit in $g_{123}$.  The trick is the same: observe that 
$g_{123}=g_{12}/\!\gcd(3,g_{12})=c_1/(\gcd(3,c_1)\gcd(2,c_1))$; from this we have
\vspace{-2pt}
\begin{align}
\gcd(3g_{123}/g_{12}, g_{123}) &= 
\gcd(3/\!\gcd(3,c_1),c_1/(\gcd(3,c_1)\gcd(2,c_1))) \nn\\
&= \gcd(3/\!\gcd(3,c_1),c_1/\!\gcd(3,c_1))=1. \nn
\end{align}
The homology of $D$ is thus, again, $\Z/g_{123}$.  The full $\Tor$ group
follows.

Despite the increased complexity of the $K$-homology rings of $\Omega E_7$
and $\Omega E_8$, the $\Tor$ calculations in these cases
are no more elaborate than for
the other exceptional groups.  The presentation in Proposition~\ref{proppres}
suggests an appropriate Tate resolution and the $\Tor$ group over
$K.\Omega E_7$ is given by the homology of the DGA
\begin{multline}
(\Z\langle T_1, \ldots, T_{10} \rangle \{S_1,S_2,S_3\};
dT_i=c_i, dS_1=(c_1 - 1)T_1 - 2T_2, \\
dS_2=(c_1 - 2)T_2 - 3T_3, 
dS_3=(c_2 + 1)T_2 - (c_2 + c_3) T_1 - 2T_4). \nn
\end{multline}
Using the $F_4$ computation, we see that the homology of the 
$\{T_1,T_2,T_3,T_4,S_1,S_2\}$ subalgebra is $\Z/g_{1234}\langle y_4
\rangle$, where $y_4$ is $-(g_{123}/g_{1234}) T_4$ modulo lower
terms.  The homology of the $\{T_1,T_2,T_3,T_4,S_1,S_2,S_3\}$
subalgebra is therefore the homology of
$$(\Z/g_{1234}\langle y_4 \rangle\{S_3\}; dS_3=(2g_{1234}/g_{123})y_4).\nn$$
We observe that $2g_{1234}/g_{123}$ is a unit in $\Z/g_{1234}$
and so the homology of this subalgebra is $\Z/g_{1234}$ concentrated
in degree zero.  The full $\Tor$ group is finally $\Z/g_{1..10}\langle 
y_5, y_6, \ldots, y_{10} \rangle$.  In this calculation it is critical
that the third relation in the presentation of 
$K.\Omega E_7$ is $b(b+1) - a(b+c) - 2d$ and 
not the expected $c(a-3)-4d$.  The latter relation would produce a
differential $dS_3=(4g_{1234}/g_{123})y_4$ and thereby (because
$4g_{1234}/g_{123}$ is not always a unit in $\Z/g_{1234}$) a plethora of
nontrivial higher torsion.

The $\Tor$ computation for $E_8$ is entirely analogous.  The Tate resolution
is dictated by the presentation in Proposition~\ref{proppres} and the
necessary combinatorial fact is that $5g_{12345}/g_{1234}$ is a unit
in $\Z/g_{12345}$.

\subsection{Proof of Theorem~\ref{thmtkh}} \label{tkh}

\vspace{-8pt}
We can now establish the bulk of our main theorem.
We assume the computation of the torsion group for $Spin(n)$, which is
carried out in section~\ref{torSOn}:
$$\Tor^{K.\Omega Spin(n)}(K.*,(K.*)_{\tau}) =
\Lambda[x_1, \ldots, x_{n-1}] \tensor \Z/c(Spin(n),k). \nn$$

Though we have treated $K$-homology as $\Z/2$-graded in our $\Tor$ 
computations,
properly it is
$\Z$-graded, and the $E^2$ term of the Rothenberg-Steenrod spectral
sequence has the appearance:
\begin{equation} \nn
\xymatrix@C=20pt@R=3pt{
\Tor_0^{K_0 \Omega G}(\Z,\Z_{\tau}) & \Tor_1^{K_0 \Omega G}(\Z,\Z_{\tau}) &
\Tor_2^{K_0 \Omega G}(\Z,\Z_{\tau}) & \cdots \\
0 & 0 & 0 & \\
\Tor_0^{K_0 \Omega G}(\Z,\Z_{\tau}) & \Tor_1^{K_0 \Omega G}(\Z,\Z_{\tau}) &
\Tor_2^{K_0 \Omega G}(\Z,\Z_{\tau}) & \cdots \\
0 & 0 & 0 & \\
\Tor_0^{K_0 \Omega G}(\Z,\Z_{\tau}) 
\save[]+<-1.6cm,-.4cm> \ar@{-}[uuuu]+<-1.6cm,.4cm>
\restore
& \Tor_1^{K_0 \Omega G}(\Z,\Z_{\tau}) &
\Tor_2^{K_0 \Omega G}(\Z,\Z_{\tau}) & \cdots
}
\end{equation}
In our cases these torsion groups are generated in $\Tor$-degree 1; the
(homological) differentials vanish on the generators and thus the
spectral sequence collapses at the $E^2$ term.

We show that there are no additive extensions.  We have established
that the $E^{\infty}$ term of the spectral sequence is
a cyclic group, say $\Z/c$, tensor an exterior algebra.  The filtration
is homological, so the subgroup $(\Z/c)\{1\} \subset E^{\infty}$ generated by 
the identity element of the torsion group $\Tor^{K.\Omega G} = E^2 = 
E^{\infty}$ is actually a subgroup of the $K$-homology $K^{\tau}(G)$.
The construction of the spectral sequence shows that this 
identity element in the torsion group corresponds to the identity element
in the $K$-homology.  The identity element $1 \in K^{\tau}(G)$ is
therefore killed by multiplication by $c$, and so the entire $K$-homology ring 
is $c$-torsion, as desired.

For degree reasons, the only possible multiplicative extension is
$y_i^2=d \in K^{\tau}(G)$; that is the square of the $K$-homology
class represented by an exterior generator $y_i \in \Tor_1$ could be a
constant integer, an element of $\Tor_0$.  However, by construction
the exterior classes $y_i$ are represented by {\em{reduced}} classes in
$K^{\tau}(G)$ and so their squares are also certainly reduced,
eliminating the possibility of multiplicative extensions.
\hfill $\Box$

\section{Generating Varieties, the Cyclic Order of $K^{\tau}G$,
and $\Tor^{K.\Omega Spin(n)}(\Z,\Z_{\tau})$} \label{gvco}

The twisted $K$-homology of a simple Lie group is an exterior algebra
tensor a cyclic group.  The order of this cyclic group depends on
the twisting class and is, as yet, determined by a mysterious set of
constants.  We will see that this cyclic order of the twisted 
$K$-homology $K^{\tau}G$ is the greatest common divisor 
of the dimensions of a particular set of representations of $G$.  
The main ingredient in computing the cyclic order is a detailed 
understanding of the twisted module structure of $\Z_{\tau}$, 
that is of the twisting map $K.\Omega G \overset{\tau(k)}{\ra} K.*$.
Bott's theory of generating varieties allows us to produce explicit
representatives of classes in $K.\Omega G$, as fundamental classes
of complex algebraic varieties, and thereby to describe the twisting
map.

\subsection{Generating Varieties and Holomorphic Induction} \label{gvhi}

\subsubsection{Bott Generating Varieties}

A generating variety for $\Omega G$ is, for us, a space $V$
and a map $i:V \ra \Omega G$ such that the images $i_*(H.V)$
and $i_*(K.V)$ of the homology and $K$-theory of $V$ generate
$H.\Omega' G$ and $K.\Omega' G$, respectively, as algebras,
where $\Omega' G$ is the identity component of $\Omega G$.
In~\cite{bottsllg}, Bott produced a beautiful, systematic family 
of generating varieties of the form $G/H$, as we now describe; 
these particular homogeneous spaces are better known as coadjoint 
orbits and as such are smooth complex algebraic varieties with 
an even-dimensional cell decomposition.

\begin{figure}[h]
\begin{center}
\caption{Generating Coweights for Rank 2 Lie Groups} \vspace{6pt}
\epsfig{figure=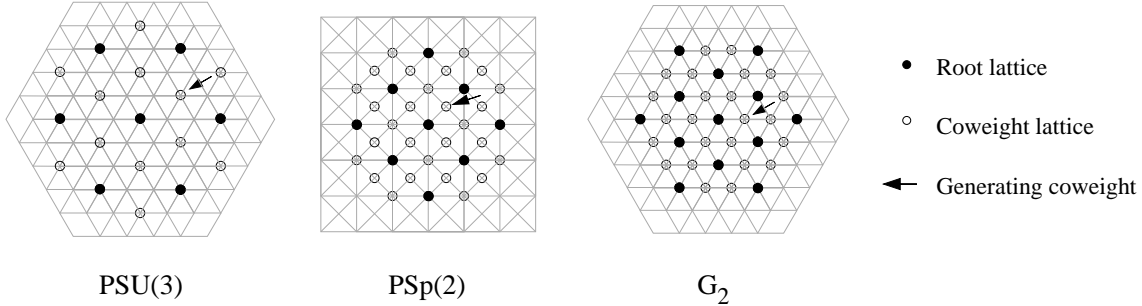,height=1.6in} \label{roots}
\end{center}
\end{figure}

We briefly review Bott's construction.  Let $G$ be a compact 
and connected but not necessarily simply connected Lie group.  
Denote by $\Gamma_G = \ker(\exp:\ttt \ra T)$ the coweight 
lattice of $G$; we do not distinguish between a coweight and the
corresponding circle in $G$.  A coweight $\ell \in \Gamma_G$ 
is called generating if for every root $r \in \ttt^*$ of $G$, 
there is an element $w$ of the Weyl group such that
$r(w \cdot \ell)=1$.  Note that the coweight lattice $\Gamma_G$
of the group is contained in the coweight lattice $\Gamma_W$ of
the Lie algebra, which is also the coweight lattice of the
adjoint form of $G$.  Though a group may not have a generating
coweight, its adjoint form always will.
The simple rank 2 groups with generating coweights, namely $PSU(3)$, 
$PSp(2)$, and $G_2$, are illustrated in Figure~\ref{roots}.

Suppose $\ell \in \Gamma_G$ is a generating coweight for $G$, and let
$C(\ell) \subset G$ denote the pointwise centralizer of the corresponding
circle; note that $C(\ell)$ can also be described as the image under the 
exponential map of the subalgebra of $\g$ generated by the root 
spaces associated to roots $r$ perpendicular to $\ell$, that is to roots 
where $r(\ell)=0$.  The map
\begin{align}
G &\ra \Omega G \nn \\
g &\mapsto g \ell g^{-1} \ell^{-1} \nn
\end{align}
descends to a map on cosets $G/C(\ell) \ra \Omega G$.  The main theorem,
which is due to Bott in homology and to Clarke~\cite{clarke} in
$K$-theory, is that $G/C(\ell)$ is a generating variety for $\Omega G$.

Suppose $G$ is simply connected and $\ell$ is a generating circle for 
$PG$, the adjoint form of $G$.  Then $PG/C_{PG}(\ell)=G/C_G(\widetilde{\ell})$
where $\widetilde{\ell}$ denotes a loop in $G$ covering $\ell$.
The composite
$$G/C_G(\widetilde{\ell}) = PG/C_{PG}(\ell) \ra \Omega' PG = \Omega G \nn $$
is therefore a generating variety for $\Omega G$.  For example, the generating
varieties corresponding to the marked coweights in Figure~\ref{roots}
are $SU(3)/U(2)$, $Sp(2)/U(2)$, and $G_2/U(2)$ respectively.  In
general there may be more than one Bott generating variety for a group;
we list a Bott generating variety for each of the classical groups in
the following table:

\begin{center}
\begin{tabular}{l|l}
Group & Generating Variety \\
\hline
$SU(n+1)$ & $SU(n+1)/U(n)$ \\
$Spin(2n+1)$ & $Spin(2n+1)/(Spin(2n-1) \times_{\Z/2} Spin(2)) $ \\
$Sp(n)$ & $Sp(n)/U(n)$ \\
$Spin(2n)$ & $Spin(2n)/(Spin(2n-2) \times_{\Z/2} Spin(2))$ \\
\end{tabular}
\end{center}
\vspace{5pt}
\nid Here the $\Z/2$ action on $Spin(n)$ is the one whose
quotient is $SO(n)$.

\vspace{5pt}

{\begin{center} $\sim\sim\sim\sim\sim$ \end{center}} \vspace{6pt}

We need to compute the twisted map $K.\Omega G \overset{\tau(k)}{\ra}
K.*$.  To this end we want to represent the algebra generators of
$K.\Omega G$ in a way that allows us to compute their twisted images.
In our computations we utilize generating varieties to represent
algebra generators in three independent ways; we refer to these
briefly as representing them via subvarieties, via an evaluation dual
basis, and via a Poincare dual basis.

In some cases we have a sufficiently explicit handle on the generating 
variety $V$ for $\Omega G$ that we can describe a collection of maps
$W_i \ra V \ra \Omega G$ such that $W_i$ is a $K$-oriented manifold
and the images in $K.\Omega G$ of the $K$-homology fundamental classes
of the $W_i$ are the desired algebra generators; frequently, though
not always, the $W_i$ are subvarieties of the generating variety $V$.
Variants of this `subvariety' representation are used for $SU(n)$,
$G_2$, and $Sp(n)$ in sections~\ref{genSUnG2} and~\ref{genSpn}.

The $K$-cohomology of the Bott generating variety $V$ is easily determined
from the representation theory of $G$.  Specifically, if the Bott generating
variety $V$ is the quotient $G/H$ with $G$ simply connected,
then $K^{\bdot}V = R[H]/i^*(I[G])$, where $i:H \ra G$ is the inclusion
and $I[G]$ is the augmentation ideal of the representation ring $R[G]$.
If there is a minor miracle and we can write down a 
clean basis for this ring, then we can take an evaluation
dual basis for the $K$-homology $K.V$; the image of this basis in
$K.\Omega G$ will generate as an algebra and the twisting map will be
easily computable.  This is the approach taken for $Spin(n)$ in section
\ref{torSOn}.

More commonly, any apparent basis for the $K$-cohomology of the generating
variety is quite haphazard.  In this case we consider the Poincare dual
basis (to some chosen basis) for the $K$-homology $K.V$.  Again the
images of these classes in $K.\Omega G$ will generate, but computing
their twisted images requires a bit more work.  Specifically, we use
holomorphic induction to write these images in terms of the dimensions
of irreducible representations of $G$, as described in detail in the
next section.  This Poincare dual approach is the one that provides
a general procedure and is the subject of section~\ref{cycgen}.

\subsubsection{The Twisting Map via Holomorphic Induction} \label{tihi}

In section~\ref{trsss} we described the $K.\Omega G$-module structure
on $\Z_{\tau}$ by the twisting map
\begin{align}
K.\Omega G &\overset{\tau(k)}{\ra} K.* \nn \\
x &\mapsto \langle L^k,x \rangle, \nn
\end{align}
where $L \in K^{\bdot}\Omega G$ is a generating line bundle.
The purpose of this section is to outline the computation of the twisted 
image $\lb L^k,x \rb$ when $x$ is represented as the image of the 
Poincare dual of a bundle on an appropriate $K$-oriented manifold.

Let $i:W \ra \Omega G$ be a map from a $K$-oriented manifold $W$ to
$\Omega G$ and let $\eta \in K^{\bdot} W$ be a bundle on $W$ such that
$i_*(\D\eta)=x \in K.\Omega G$; here, $\D$ denotes the Poincare duality 
map.  We first translate the evaluation $\lb L^k,x \rb$ into a pushforward on
$W$:
$$ \lb L^k, i_*(\D\eta) \rb = \lb i^*(L^k), \D\eta \rb = \lb i^*(L^k), 
\eta \cap [W] \rb
= \lb i^*(L^k) \cup \eta, [W] \rb. \nn $$
The third equality is simply
$$ \lb i^*(L^k), \eta \cap [W] \rb = \pi_*(i^*(L^k) \cap (\eta \cap [W])) = 
\pi_*((i^*(L^k) \cup \eta) \cap [W]) = \lb i^*(L^k) \cup \eta, [W] \rb, \nn $$
where $\pi:W \ra *$ denotes projection.  We are thereby reduced to
computing $K$-theory pushforwards $\lb \mu,[W] \rb =\pi_!(\mu)$.

Suppose that, as in our computations it always is, $W$ is a homogeneous
space $G/H$ which is a K\"ahler manifold, and $\mu \in K^{\bdot}(G/H)$ is the bundle
associated to an irreducible representation of $H$.  (If our original bundle $\mu$ is
not irreducible, we decompose it into irreducible components and work
with each component separately.)  In this case $\mu$ lifts to a $G$-equivariant
bundle, also denoted $\mu$, and the pushforward $\pi_!(\mu)$ is
the dimension of the equivariant pushforward $\pi_!^G(\mu)$.  This equivariant
pushforward $\pi_!^G: R[H] \cong K_G^{\bdot}(G/H) \ra K_G^{\bdot}(\ast) \cong R[G]$
is usually referred to as holomorphic induction.  Thanks to the
Atiyah-Segal fixed point theory~\cite{as}, pushforwards in equivariant $K$-theory,
particularly those involving homogeneous spaces, are easily computable.
To avoid reviewing the whole of Atiyah-Segal's theory, we describe only the
form it takes in the context of holomorphic induction.

Let $G$ be a compact connected simply connected Lie group and $H$ a 
centralizer of a circle in $G$; in particular $H$ and $G$ share a maximal 
torus and their weight lattices coincide.  In this situation, as remarked
earlier, the $K$-theory of the quotient $G/H$ is simply the quotient of
representation rings: $K^{\bdot}(G/H) = R[H]/i^*(I[G]).$  We may further
assume that we have chosen an order on the roots of $G$ such that
$H$ is generated by a subset of the simple roots of $G$; in particular
this determines standard Weyl chambers for $G$ and $H$.  (That
there is such a choice of order is the content
of Wang's theorem and depends on $H$ being the centralizer of a torus
in $G$; see Bott~\cite{botthvb}.)  Let $\mu$ denote 
simultaneously a weight in the Weyl chamber of $H$, the corresponding
irreducible representation of $H$, and the associated bundle on $G/H$.
Let $\rho$ denote half the sum of the positive roots of $G$, and let
$S$ denote the union of the hyperplanes perpendicular to the roots of $G$.
Further, for a weight $\omega$ of $G$, let $T(\omega)$ denote the unique 
weight in the Weyl chamber of $G$ that is the image of $\omega$ under an
element of the Weyl group.  The index $\ind(\omega)$ of a weight $\omega$ 
not in $S$ is the number of hyperplanes of $S$ intersecting a straight line
connecting $\omega$ and $T(\omega)$.  Holomorphic induction on the
representation $\mu$ is described by the following well known theorem:
\vspace{4pt}
\begin{theorem} \label{botthi}
In the above situation,
\begin{enumerate}
\item The character of the representation $\pi_!^G(\mu) \in R[G]$ is
$\frac{\sum\limits_{w \in W} \det(w) (\mu \cdot \rho)^w}{\sum\limits_{w \in W} \det(w) (\rho)^w}$.
Here $\mu \cdot \rho$ is the one-dimensional representation with weight $\mu+\rho$, and
$W$ is the Weyl group of $G$.
\item If $\mu+\rho \notin S$ then $\pi_!^G(\mu) \in R[G]$ is $(-1)^{\ind(\mu+\rho)}$ times
the irreducible representation of $G$ with highest weight $T(\mu+\rho)-\rho$.
\end{enumerate}
\end{theorem}
\vspace{4pt}
\nid The first part of the theorem is a consequence of the Atiyah-Segal fixed point formula~\cite{as},
and the second part follows from the first using the Weyl character formula.  Strictly speaking, we
only need the first part, but it will be convenient, using the second part of the theorem, to refer to appropriate pushforwards as irreducible representations.

The nonequivariant pushforwards follow immediately.  When $\mu + \rho \in S$ we say that $\mu$ is singular.  Thus, when $\mu$ is singular, $\pi_!(\mu) = 0$, and otherwise
$$\pi_!(\mu) = (-1)^{\ind(\mu+\rho)} \dim([T(\mu+\rho)-\rho]_G), \nn$$
where $[-]_G$ denotes the irreducible representation of $G$ with the specified highest weight.
The character of that representation is given in~\ref{botthi} above, and in any particular case its dimension is easily computed using the Weyl dimension formula.  This method provides a systematic approach to computing the twisting map on a class represented by the image of the Poincare dual to a $K$-cohomology class of an appropriate homogeneous space.  We proceed to specific examples.

\vspace{-6pt}
\subsection{Subvarieties of $\Omega SU(n+1)$ and $\Omega G_2$} \label{genSUnG2}

\vspace{-10pt}
A generating variety for $\Omega SU(n+1)$ is $SU(n+1)/U(n)=\CP^n
\overset{i}{\ra} \Omega SU(n+1)$, and the induced map in homology is
\vspace{-4pt}
\begin{align}
\widetilde{H}.\CP^n = \Z\{z_1, \ldots, z_n\} &\ra \Z[w_1, \ldots, w_n] =
H.\Omega SU(n+1) \nn \\
z_i &\mapsto w_i. \nn
\end{align} \vspace{-20pt}

\nid Here the classes $z_i$ are represented by the fundamental homology
classes of the subvarieties $\CP^i \subset \CP^n$.  The
Atiyah-Hirzebruch spectral sequence collapses for both $\Omega SU(n+1)$
and $\CP^n$ and there are no extensions.  In particular
$K.\Omega SU(n+1)$ is polynomial on $n$ generators, as previously
noted, and $\widetilde{K}.\CP^n$ is free abelian of rank n.
\vspace{8pt}
\begin{lemma} \label{SUnlemma}
The set $\{[\CP^i]\}_{i=1}^n$ of fundamental $K$-homology classes of
the subvarieties $\CP^i \subset \CP^n$ forms a basis for 
$\widetilde{K}.\CP^n$.
\end{lemma} \pagebreak
\begin{proof}
By induction it is enough to show that under the projection
$K.\CP^i \ra K.(\CP^i,\CP^{i-1}) = \Z$
the fundamental class of $\CP^i$ maps to a generator.  This
follows immediately from the naturality of Poincare duality,
\vspace{-6pt}
\begin{equation} \nn
\xymatrix@-5pt{
K.\CP^i \ar[r] \ar@{=}[d] & K.(\CP^i, \CP^{i-1}) \ar@{=}[d] \\
K^{\bdot}\CP^i \ar[r] & K^{\bdot}(\CP^i - \CP^{i-1}),
}
\end{equation}

\vspace{-8pt}
\nid because the unit in $K^{\bdot}\CP^i$ certainly maps to a generator of
$K^{\bdot}(\CP^i - \CP^{i-1})$.
\end{proof}
\nid The images $i_*([\CP^i])$ generate $K.\Omega SU(n+1)$ as an
algebra and we may therefore take $\{x_i=i_*([\CP^i]) - 1\}$ to be
the reduced polynomial generators.

We now have to evaluate the pushforward $\lb L^k,[\CP^i] \rb$, where
we use $L$ to denote the generating line bundle on $\CP^i=
SU(i+1)/U(i)$; this $L$ is of course the pullback of the generating 
line bundle on $\Omega SU(n+1)$.  The bundle $L$ corresponds to
an irreducible representation of $U(i)$, thus to a weight of $U(i)$
and so a weight, also denoted $L$, of $SU(i+1)$; this weight $L$ is
in the Weyl chamber of $SU(i+1)$.  The irreducible representation of 
$SU(i+1)$ corresponding to $L$ is the dual of the standard 
representation.  (That it is the dual of the standard
representation and not the standard representation is the effect of
a sign choice---see the remark at the end of this section).  It happens 
that the $k$-fold symmetric power of this representation is irreducible, 
and so the dimension of the irreducible representation corresponding 
to $L^k$ is $\binom{k+i}{i}$.  The image of $x_i = i_*([\CP^i]) - 1 
\in K.\Omega SU(n+1)$ in $\Z_{\tau}$ is therefore $c_i=\binom{k+i}{i} - 1$, 
and the cyclic order of $K^{\tau}(SU(n+1))$ is
$$c(SU(n+1),k) = \gcd\left\{\binom{k+1}{1} - 1, \binom{k+2}{2} - 1,
\ldots, \binom{k+n}{n} - 1\right\}.$$

{\begin{center} $\sim\sim\sim\sim\sim$ \end{center}} \vspace{6pt}

The procedure for calculating the cyclic order of $K^{\tau} G_2$
is similar: we find fundamental class representatives for algebra generators 
of the homology of $\Omega G_2$ and then show that the corresponding 
$K$-homology fundamental classes also generate.
The map $\Omega G_2 \ra \CP^{\infty}$ classifying
the generating line bundle is a homology equivalence through degree 4.
Using this, the Serre spectral sequence for $\Omega SU(3) \ra \Omega G_2
\ra \Omega S^6$ shows that 
$$H.\Omega G_2 \cong \Z[a_2,a_4,a_{10}]/a_2^2=2a_4.\nn$$
The composition
$$\CP^2 \ra \Omega SU(3) \ra \Omega G_2 \ra \CP^{\infty}$$
is simply the inclusion and as such, $a_2$ and $a_4$ in $H.\Omega G_2$
are represented respectively by the fundamental classes $[\CP^1]$ and 
$[\CP^2]$.  The Bott generating variety for $\Omega G_2$ is
$G_2/U(2)$, where the $U(2)$ in
question is included in $G_2$ along a pair of complex-conjugate
short roots.  The manifold
$G_2/U(2)$ has dimension 10 and the image of its homology generates
$H.\Omega G_2$; we may therefore choose $a_{10}$ to be the image of
the fundamental homology class $[G_2/U(2)]$.

The Atiyah-Hirzebruch spectral sequence for $\Omega G_2$ collapses,
and the low-degree equivalence between $\Omega G_2$ and $\CP^{\infty}$
resolves the extension.  The $K$-homology $K.\Omega G_2$ is thereby
isomorphic to $\Z[a,b,x_3]/(a^2+3a=2b)$.
\vspace{8pt}
\begin{lemma}
The reduced algebra generators of $K.\Omega G_2 \cong \Z[a,b,x_3]/(a^2+3a=2b)$
may be taken to be the reduced fundamental $K$-homology classes
$[\CP^1]-1$, $[\CP^2]-1$, and $[G_2/U(2)]-1$ respectively.
\end{lemma}
\begin{proof}
Let $(G_2/U(2))_8$ denote the 8-skeleton of the generating variety,
that is everything except the top cell.  As in Lemma~\ref{SUnlemma},
the fundamental $K$-homology class of $G_2/U(2)$ maps to a generator
of $K.(G_2/U(2),(G_2/U(2))_8)$.  Comparing the Atiyah-Hirzebruch
spectral sequences for $G_2/U(2)$ and $\Omega G_2$ we see that
$[G_2/U(2)]$ lives in filtration 10 in $K.\Omega G_2$ and projects
to the generator $a_{10}$ in $H_{10} \Omega G_2$.  The fundamental
$K$-homology classes $[\CP^1]$ and $[\CP^2]$ certainly project
to the generators $a_2$ and $a_4$ respectively in the appropriate
filtration quotients, and this completes the proof.
\end{proof}
\nid We remark that these algebra generators differ by a change of
basis from those implicitly chosen in section~\ref{torG2F4E6}
and this explains the difference in the relation; the $\Tor$
computation and the cyclic order are not affected by the change.

We need only compute the pushforward $\lb L^k,[G_2/U(2)] \rb$.  The
bundle $L$ corresponds to the shortest weight $\mu$ perpendicular to
the roots of $U(2)$; as a weight of $G_2$, $\mu$ is the long root
of $G_2$ in the Weyl chamber.  The pushforward is therefore the
dimension of the irreducible representation of $G_2$ with highest
weight $k \mu$.  By the Weyl dimension formula (see for 
example~\cite{fh}) this is
$$\dim([k \mu]_{G_2}) = \frac{(k+1)(k+2)(2k+3)(3k+4)(3k+5)}{120}. \nn$$
The cyclic order of $K^{\tau} G_2$ is finally
$$c(G_2,k) = \gcd\left\{\binom{k+1}{1} - 1, \binom{k+2}{2} - 1,
\frac{(k+1)(k+2)(2k+3)(3k+4)(3k+5)}{120} - 1\right\}. \nn$$

{\begin{center} $\sim\sim\sim\sim\sim$ \end{center}} \vspace{6pt}

A remark on signs is in order.  If we have chosen a generating line bundle
$L$ on $\Omega G_2$ a priori, the weight corresponding to $L$ may be $-\mu$
instead of $\mu$ as claimed above.  The dimension resulting from
holomorphic induction on the weight $k (-\mu)$ is wildly different from
the dimension associated to $k (\mu)$, and this might be cause for
worry.  However, the greatest common divisor is in all cases unaffected
by the change.  The easiest way around this ambiguity is to chose $L$
such that the corresponding weight is $\mu$ and not $-\mu$; we must then
pick the generating variety $\CP^2$ for $\Omega SU(3)$ in such a way
that the given $L$ corresponds there to the dual of the standard
representation (and not to the standard representation) as described in
the discussion of $SU(n+1)$ above---this is easily accomplished.  Similar
remarks apply to all our computations and we make convenient sign choices
without comment.

\subsection{Generating Varieties for $\Omega Sp(n)$} \label{genSpn}

The homology and $K$-homology of $\Omega Sp(n)$ are polynomial in $n$
generators.  The natural Bott generating variety for $\Omega Sp(n)$ 
is $Sp(n)/U(n)$, which has homology and $K$-homology of rank $n^2+n$.
Identifying the $n$ elements which generate therefore requires more
doing---we return to this question later.  Luckily, $\Omega Sp(n)$
has smaller generating varieties---see~\cite{hopkins,mitchell}; 
in particular we work with $(\CP^{2n-2})^{L^2}$, the Thom complex of 
the square of the tautological bundle.

Let $P_i(V)$ or $P(V)$ denote the projectivization of the bundle $V$
on $\CP^i$; note that we can rewrite our generating variety
$V(n)=(\CP^{2n-2})^{L^2}$ as $P_{2n-2}(L^2+1)/P_{2n-2}(L^2)$.
We think of the quotient map $P(L^2+1) \ra P(L^2+1)/P(L^2)$ as a
resolution of our (quite singular) generating variety, and
we represent homology and $K$-homology classes in $V(n)$ (and
thus in $\Omega Sp(n)$) as the images of fundamental classes
of subvarieties of $P(L^2+1)$.  The reduced homology of
$V(n)$ is free of rank one in each even degree between $2$
and $4n-2$, and the degree $2i$ group is generated by the
image of the fundamental class $[P_{i-1}(L^2+1)]$.  In particular,
the algebra generators $\{a_{4i-2}\}$ of $H.\Omega Sp(n) = 
\Z[a_2, a_6, a_{10}, \ldots, a_{4n-2}]$ are represented by
the fundamental classes $[P_{2(i-1)}(L^2+1)]$, for $1 \leq i \leq n$.

The $K$-homology situation is the same.
\vspace{8pt}
\begin{lemma}
The reduced polynomial generators of the $K$-homology $K.\Omega Sp(n)
\cong \Z[x_1, \ldots, x_n]$ can be taken to be the reduced
$K$-homology fundamental classes $f_*[P_{2(i-1)}(L^2+1)] - 1$,
$1 \leq i \leq n$; here $f$ is the composite
$$P_{2(i-1)}(L^2+1) \ra P_{2(n-2)}(L^2+1) \ra 
P_{2(n-2)}(L^2+1)/P_{2(n-2)}(L^2) \ra \Omega Sp(n). \nn$$
\end{lemma}
\vspace{5pt}
\nid The $K$-homology fundamental classes map, in the appropriate
filtration quotients, to the homology fundamental classes; the
proof is the same as for $SU(n+1)$ and $G_2$.  

To evaluate
the twisting map, specifically to calculate $\lb L^k,f_*[P_{2(i-1)}
(L^2+1)] \rb$,
we need to identify the bundle $f^*(L^k)$.  We do this by writing
down a bundle on $P(L^2+1)=P_{2(i-1)}(L^2+1)$ that is trivial on 
$P(L^2)=P_{2(i-1)}(L^2)$, and show that the corresponding bundle 
on the quotient $V(i)$ is the pullback $f'^*(L)$ where $f'$ is the 
inclusion $V(i) \ra V(n) \ra \Omega Sp(n)$.  Let $\gamma$ be the 
tautological bundle on the total space $P(L^2+1)$ and let 
$\pi:P(L^2+1) \ra \CP^{2(i-1)}$ be the bundle projection.  The subspace 
$P(L^2)$ is of course just the base $\CP^{2(i-1)}$ and so $\gamma$ 
restricts to $\pi^*(L^2)|_{P(L^2)}$ on $P(L^2)$.  In particular 
then, the bundle $\gamma \tensor \pi^*(L^{-2})$ is trivial on the subspace
$P(L^2)$ and so pulls back from a bundle $\phi$ on $V(i)$.  To see that
$\phi$ is equal to $f'^*(L)$ (up to our usual sign ambiguity), and therefore 
that $\gamma \tensor \pi^*(L^{-2})=f^*(L)$, it is enough to check that the 
first Chern class of $\gamma \tensor \pi^*(L^{-2})$ is a module generator
of $H^2(P(L^2+1))=\Z\{c_1(\gamma),\pi^*(c_1(L))\}$; this much is clear.

We now compute the pushforward 
$$\lb L^k, f_*[P_{2(i-1)}(L^2+1)] \rb =
\lb (\gamma \tensor \pi^*(L^{-2}))^k,
[P_{2(i-1)}(L^2+1)] \rb \nn.$$  First pushforward along the fibres:
$$\lb \gamma^k \tensor \pi^*(L^{-2k}), [P_{2(i-1)}(L^2+1)] \rb
=\lb \Sym^k(L^2+1) \tensor L^{-2k}, [\CP^{2(i-1)}] \rb. \nn$$
This is a parameterized version of the pushforward 
$$\lb \gamma_{\rm taut}^k,[\CP^n] \rb= \lb \gamma_{\rm taut}^k, [P(\C^{n+1})]
\rb =\dim(\Sym^k(\C^{n+1}))=\binom{k+n}{k} \nn$$
used in the preceding section.  Next 
$$\Sym^k(L^2+1) \tensor L^{-2k}
=L^{-2k} + L^{-2k+2} + \ldots + 1 \nn$$ and so
$$\lb \Sym^k(L^2+1) \tensor L^{-2k}, [\CP^{2(i-1)}] \rb =
\binom{-2k+2(i-1)}{2(i-1)} + \binom{-2k+2+2(i-1)}{2(i-1)} + \ldots + 1. \nn$$
Here we use $\binom{a+b}{b}$ to denote 
$\frac{(a+b)(a+b-1)\ldots(a)}{(b)(b-1)\ldots(1)}$ even when $a+b$ is negative;
we implicitly observe that this expression does give the correct pushforward
even when the bundle, as in the case of $L^{-2k+2l}$, corresponds to a weight
that is not in the Weyl chamber of $SU(2(i-1)+1)$.  This finishes our
calculation of the cyclic order of $K^{\tau}Sp(n)$:
$$c(Sp(n),k) = \gcd\left\{
\sum\limits_{-k \leq j \leq -1} \binom{2j+2(i-1)}{2(i-1)} : 1\leq i \leq n
\right\}. \nn$$

{\begin{center} $\sim\sim\sim\sim\sim$ \end{center}} \vspace{6pt}

It would be more natural to express the cyclic order of
$K^{\tau}Sp(n)$ in terms of the dimensions of irreducible
representations of symplectic groups.  This is possible if we work with
subvarieties of the Bott generating variety $V=Sp(n)/U(n)$.  There
is a natural collection of $n$ subvarieties of $V$, namely
$\{Sp(i)/U(i)\}$.  It is not the case that the fundamental homology
classes of these subvarieties represent algebra generators for
$H.\Omega Sp(n)$; indeed, the algebra generators are in dimensions
$\{4i-2\}$, while these subvarieties have dimensions $\{i^2+2\}$.
It is therefore remarkable that the $K$-homology fundamental classes
of these subvarieties do appear to generate the $K$-homology of
$\Omega Sp(n)$.
\vspace{8pt}
\begin{conject}
The $K$-homology ring $K.\Omega Sp(n)$ is polynomial on the
classes represented by the reduced $K$-homology fundamental
classes $[Sp(i)/U(i)]-1$, for $1 \leq i \leq n$.
\end{conject}
\vspace{10pt}
\nid Using the Weyl character formula, this immediately gives a description
of the cyclic order:
$$c(Sp(n),k) = \gcd{\textstyle \left\{ \frac{\left(\prod\limits_{1 \leq j <
l \leq i} (l-j)(2k+2i+2-(j+l))\right)\left(\prod\limits_{1 \leq j \leq i} 
(k+i+1-j)\right)}
{(2i-1)! (2i-3)! \ldots 3! 1!} : 1\leq i \leq n\right\}.} \nn $$
These $\gcd$'s agree with those determined using the generating
variety $(\CP^{2n-2})^{L^2}$.

\subsection{The $\Tor$ Calculation for $Spin(n)$} \label{torSOn}

We now pay our debt to the proof of Theorem~\ref{thmtkh} by calculating
$\Tor^{K.\Omega Spin(n)}(\Z,\Z_{\tau})$; in the process we determine
the cyclic order of $K^{\tau} Spin(n)$.  For the other simple groups, we 
were able to calculate the $\Tor$ group
without knowing the map $K.\Omega G \ra (K.*)_{\tau}$ and we determined,
after the fact, the structure of this twisting map.  The ring
$K.\Omega Spin(n)$ is too complicated to permit this a-priori $\Tor$
calculation; we must first identify algebra generators of
$K.\Omega Spin(n)$ and compute the twisting map.  It happens that
the reduced $K$-cohomology of the Bott generating variety for $Spin(n)$
admits a particularly simple 
representation-theoretic basis,
and an evaluation dual basis maps to a set of algebra generators in
$K.\Omega Spin(n)$.  Once we know the twisted pushforwards of these
algebra generators, the $\Tor$ computation becomes tractable.

We concentrate on the odd $Spin$ groups; at the end of the section we delineate
the corresponding steps for the even $Spin$ groups.
The structure of the $K$-homology ring of $\Omega Spin(2n+1)$ was described
by Clarke~\cite{clarke}:
\begin{align}
K.\Omega Spin(2n+1) &= \frac{\Z[\sigma_1,\sigma_2,\ldots,\sigma_{n-1},
2\sigma_{n},2\sigma_{n+1}+\sigma_{n},\ldots,2\sigma_{2n-1}+\sigma_{2n-2}]}
{(\rho_1,\ldots,\rho_{n-1})}, \nn \\
\rho_k &= \sigma_k^2 + \sum\limits_{i=0}^{k-1} (-1)^{k-i}\sigma_i
\sum\limits_{j=k}^{2k-i-1}\binom{k-i-1}{j-k} (2\sigma_{j+1} + \sigma_j). \nn
\end{align}
One can see why the a-priori $\Tor$ calculation is unlikely to be
fruitful.  The $K$-cohomology of the Bott generating variety $V=Spin(2n+1)/
(Spin(2n-1) \times_{\Z/2} Spin(2))$  is simply
the quotient of the representation ring of $Spin(2n-1) \times_{\Z/2} Spin(2)$ 
by the image of the augmentation ideal of the representation ring of 
$Spin(2n+1)$.  Clarke writes this quotient in a convenient form:
$$K^{\bdot} V = \Z[\mu,\gamma]/(\mu^n-2\gamma-\mu\gamma,\gamma^2); \nn$$
here $\mu = L - 1$ where $L$ is the generating line bundle whose $k$-th
power determines the twisting.  Note that $\mu^{2n}=0$ in this ring, and
so $(\mu, \mu^2, \ldots, \mu^{2n-1})$ is a basis for $\widetilde{K}^{\bdot}(V) \tensor
\Q$.  Letting $(\sigma_1^{'},\ldots,\sigma_{2n-1}^{'})$ be the evaluation
dual basis of $\widetilde{K}.(V)\tensor \Q$, we see that
$$ (\sigma_1^{'},\sigma_2^{'},\ldots,\sigma_{n-1}^{'},2\sigma_n^{'},
2\sigma_{n+1}^{'}
+ \sigma_n^{'}, \ldots, 2\sigma_{2n-1}^{'} + \sigma_{2n-2}^{'}) \nn$$
is a basis for $\widetilde{K}.V$; these elements map, respectively, to the given
generators of \mbox{$K.\Omega Spin(2n+1)$}.  The twisting map $K.\Omega Spin(2n+1)
\ra (K.*)_{\tau}$ takes a generator $g$ to $\lb L^k,g \rb\in\Z$.  Because
$\mu^{2n}=0$, we have
$$\lb L^k,\sigma_i^{'} \rb = \lb (\mu+1)^k,\sigma_i^{'} \rb
=\binom{k}{i}, \nn $$
and the images of our integral generators are respectively
$$\left(\binom{k}{1}, \binom{k}{2},\ldots,\binom{k}{n-1}, 2\binom{k}{n},
2\binom{k}{n+1}+\binom{k}{n},\ldots, 2\binom{k}{2n-1}+\binom{k}{2n-2}
\right).\nn$$

We can now prove that $\Tor^{K.\Omega Spin(2n+1)}(\Z,\Z_{\tau})$ is
an exterior algebra on $n-1$ generators tensor a cyclic group.
We first rewrite the above presentation of $K.\Omega Spin(2n+1)$ in a
way that suggests a propitious choice of Tate resolution.  Let 
$(a_1, \ldots, a_{2n-1})$ denote the given generators of $K.\Omega
Spin(2n+1)$.  For $i$ sufficiently large, the relation $\rho_i$
expresses the generator $a_{2i}$ in lower terms; in particular
\begin{align}
K.\Omega Spin(4n-1) &= \frac{\Z[a_1,a_2,\ldots, a_{2n-2}, a_{2n-1}, a_{2n+1},
a_{2n+3},\ldots,a_{4n-5},a_{4n-3}]}{(\rho_1,\rho_2,\ldots,\rho_{n-1})} \nn \\
K.\Omega Spin(4n+1) &= \frac{\Z[a_1,a_2,\ldots, a_{2n-2}, a_{2n-1}, a_{2n+1},
a_{2n+3},\ldots,a_{4n-3},a_{4n-1}]}{(\rho_1,\rho_2,\ldots,\rho_{n-1})}. \nn
\end{align}
The remaining relations can be written
$$\rho_i = 2 a_{2i} + r_i a_{2i-1} + \ldots, \nn$$
with $r_i$ odd and all unspecified monomials containing some $a_j$ with
$j<2i-1$, except for $\rho_1$ which is $2a_2+a_1-a_1^2$. If we can show that 
$\Tor$ over the subring $R_n=\Z[a_1,\ldots,a_{2n-2}]/(\rho_1,\ldots,
\rho_{n-1})$ is exterior on $n-2$ generators, the desired result follows.
Rather than presenting the general induction immediately, we discuss the 
first few cases explicitly.  

The case $n=1$ corresponding to $Spin(3)$ requires no
comment.  The ring $K.\Omega Spin(7)$ is $\Z[a_1,a_2,a_3,a_5]/
(2a_2+a_1-a_1^2)$.  This is reminiscent of $K.\Omega G_2$ and indeed
the Atiyah-Hirzebruch spectral sequence for the fibration
$\Omega G_2 \ra \Omega Spin(7) \ra \Omega S^7$ collapses; there are
no possible multiplicative extensions and so this confirms that
$K.\Omega Spin(7)$ is $K.\Omega G_2$ adjoin a generator in degree 6.
As in section~\ref{torG2F4E6}, the $\Tor$ group in question is
$$\Tor^{R_2}(\Z,\Z_{\tau}) = \Tor^{\Z[a_1,a_2]/(2a_2+a_1-a_1^2)}(\Z,\Z_{\tau}) 
= \Z/g_{12}. \nn$$
(Note that the generator $a_i$ of the subring $R_n$ has image under
the twisting map $c_i=\binom{k}{i}$ and as before we abbreviate
$\gcd\{c_1,c_2,\ldots,c_i\}$ by $g_{1..i}$.)

The relevant subring of $K.\Omega Spin(11)$ is
$$R_3=\Z[a_1,a_2,a_3,a_4]/(\rho_1, 2a_4+3a_3+(a_2+1)a_2+(-2a_3-a_2)a_1). \nn$$
This presentation suggests the Tate resolution
\begin{multline} \nn
\Tor^{R_3} = H(\Z \lb T_1,T_2,T_3,T_4 \rb \{S_1,S_2\}; \\
dT_i=c_i,dS_1=2T_2+(1-c_1)T_1, dS_2=
2T_4+3T_3+(c_2+1)T_2+(-2c_3-c_2)T_1). \nn
\end{multline}
The $E_1$ term of the spectral sequence associated to the filtration of this 
complex by $S_2$ is
\begin{equation} \nn
\xymatrix@C=20pt@R=3pt{
& \Z/g & \Z/g\oplus\Z/g & \Z/g \\
\Z/g & \Z/g\oplus\Z/g & \Z/g & \\
\Z/g\oplus\Z/g & \Z/g & & \\
\Z/g 
\save[]+<-1.1cm,-.4cm> \ar@{-}[uuu]+<-1.1cm,.4cm>
\ar@{-}[rrr]+<.5cm,-.4cm>
\restore
& & &
}
\end{equation}
where $g=g_{1234}$ and the generator in degree $(1,1)$ is $S_2$.  At first
blush the generators in degree $(0,1)$ have the form $t_3=(g_{12}/g_{123}) T_3
+ O(2)$ and $t_4=(g_{123}/g_{1234}) T_4 + O(3)$, where the omitted terms 
contain only terms involving ($T_2$ and $T_1$) and ($T_3$, $T_2$, and $T_1$)
respectively.  In order to determine the differential on $S_2$ we need control 
over the $T_3$ term in the generator $t_4$.  The basic observation is that if 
there exists a cocycle $t_4^{'}$ of the form $(g_{123}/g_{1234}) T_4 + O(2)$, 
then some linear combination $t_4+ct_3$ is cohomologous to $t_4^{'}$ and so we 
may take the generators in degree $(0,1)$ to be $t_4^{'}$ and $t_3$.  The 
existence
of this cocycle is ensured by the fact that $(g_{123}/g_{1234}) g_4$ is
divisible by $g_{12}$, as is easily checked.  The differential on $S_2$ is
therefore $(2g_{1234}/g_{123}) t_4^{'} + (3g_{123}/g_{12}) t_3$.  Because
the greatest common divisor of $2g_{1234}/g_{123}$ and $3g_{123}/g_{12}$ is
always 1, the torsion group is finally
$$\Tor^{R_3} = \Z/g_{1234} \lb x_4 \rb;$$
here we can choose the generator $x_4$ to be $(g_{123}/g_{1234}) T_4 + O(3)$.

The case of $Spin(15)$ proceeds similarly.  The relevant subring of
$K.\Omega Spin(15)$ is
$$R_4=\Z[a_1,a_2,a_3,a_4,a_5,a_6]/(\rho_1,\rho_2, 2a_6+5a_5+4a_4+O(3)),\nn$$
and we take the corresponding Tate resolution.  Filtering by $S_3$ we
have the spectral sequence (here condensed)
\begin{equation} \nn
\xymatrix@C=20pt@R=3pt{
\Z/g \lb x_4,x_5,x_6 \rb & \Z/g \lb x_4,x_5,x_6 \rb & \Z/g \lb x_4,x_5,x_6 \rb 
& \ldots
}
\end{equation}
The torsion $g$ is $g_{1..6}$ and the generators in degree $(0,1)$ are
\begin{align}
x_4 &= (g_{123}/g_{1234}) T_4 + O(3) \nn \\
x_5 &= (g_{1234}/g_{1..5}) T_5 + O(4) \nn \\
x_6 &= (g_{1..5}/g_{1..6}) T_6 + O(5). \nn
\end{align}
It happens that $(g_{1234}/g_{1..5})g_5$ and 
$(g_{1..5}/g_{1..6})g_6$
are both divisible by $g_{123}$; we can therefore adjust our generators so
that they are
\begin{align}
x_4 &= (g_{123}/g_{1234}) T_4 + O(3) \nn \\
x_5 &= (g_{1234}/g_{1..5}) T_5 + O(3) \nn \\
x_6 &= (g_{1..5}/g_{1..6}) T_6 + O(3). \nn
\end{align}
The differential on $S_3$ is thus $(2g_{1..6}/g_{1..5})x_6
+ (5g_{1..5}/g_{1234})x_5 + (4g_{1234}/g_{123})x_4$.  Because
$2g_{1..6}/g_{1..5}$ and $5g_{1..5}/g_{1234}$ are relatively
prime, there exist constants $z_1$ and $z_2$ so that if we set \vspace{-9pt}
\begin{align}
& y_6 = x_6 + z_1 x_4 = g_{1..5}/g_{1..6} T_6 + O(4) \nn \\
& y_5 = x_5 + z_2 x_4 = g_{1234}/g_{1..5} T_5 + O(4), \nn
\end{align}
then $\{dS_3,y_6,y_5\}$
forms a basis for the degree $(0,1)$ group.  Finally, then, the $\Tor$
group is
$$\Tor^{R_4} = \Z/g_{1..6} \lb y_5,y_6 \rb$$
as desired. \newpage

The general case is now clear.  Suppose we know that 
$$\Tor^{R_n}=\Z/g_{1..(2n-2)} \lb x_{n+1},\ldots,x_{2n-2} \rb,$$
where $x_i=(g_{1..(i-1)}/g_{1..i})T_i + O(i-1)$.
The ring $R_{n+1}$ has two additional generators $a_{2n-1}$ and
$a_{2n}$ and one additional relation $\rho_n$.  Filter the
appropriate Tate resolution by powers of $S_n$, then adjust the
generators of the degree $(0,1)$ group in the spectral sequence
so that the single generator $x_{2n}$ involving $T_{2n}$ does not
contain any terms involving $T_{2n-1}$.  This is possible because
$g_{1..(2n-2)}$ divides $(g_{1..(2n-1)}/g_{1..(2n)})
g_{2n}$.  The differential of $S_n$ then has the form
$$dS_n=(2g_{1..(2n)}/g_{1..(2n-1)})x_{2n} + (rg_{1..(2n-1)}/
g_{1..(2n-2)})x_{2n-1} + \ldots.\nn$$  As those two leading terms
are relatively prime, this ensures that $\Tor^{R_{n+1}}$ again has
the desired form.  Note that in theory there could be multiplicative
extensions in the filtration spectral sequence calculating the
$\Tor$ group, but the above procedure gives a sufficiently explicit
handle on the generating classes as to eliminate this possibility.

This completes the proof of Theorem~\ref{thmtkh} for the odd $Spin$
groups and also establishes the odd $Spin$ cyclic orders given in 
Theorem~\ref{thmcyc}.  The calculation for the even $Spin$ groups
is analogous and proceeds as follows.  The relevant $K$-homology
ring, initially described by Clarke~\cite{clarke},  is
\begin{equation}\nn
K.\Omega Spin(2n+2) = \frac{\Z[\sigma_1,\ldots,\sigma_{n-1}, \sigma_n + \epsilon, - 2\epsilon, 2\sigma_{n+1} - \epsilon, 2\sigma_{n+2} + \sigma_{n+1},\ldots,2\sigma_{2n} + \sigma_{2n-1}]}
{(\rho_1,\ldots,\rho_{n-1}, \rho_n - \epsilon^2)},
\end{equation}
where the polynomial expressions $\rho_k$ are as in the odd orthogonal case.  The $K$-cohomology
of the corresponding Bott generating variety $V_{2n+2} = Spin(2n+2)/Spin(2n) \times_{\Z/2} Spin(2)$ is
\begin{equation} \nn
K^{\bdot} V_{2n+2} = \begin{cases} 
\Z[\mu,\gamma] / (\mu^{n+1} - 2\mu\gamma-\mu^2\gamma, \gamma^2 - \mu^n \gamma, \mu^{n+1}\gamma) &\text{if $n$ is even},\\
\Z[\mu,\gamma] / (\mu^{n+1}-2\mu\gamma-\mu^2\gamma,\gamma^2)
&\text{if $n$ is odd}.
\end{cases}
\end{equation}
This presentation~\cite{clarkelet} is a slight correction of the one given in Clarke's paper.  Note that $\mu^{2n+1}=0$ in either ring, and thus 
$(1,\mu,\ldots,\mu^{n-1},\mu^n,\beta,\mu^{n+1}, \mu^{n+2}, \ldots,\mu^{2n})$ is a basis for $K^{\bdot}(V) \otimes \Q$, where $\beta = \mu^n - 2\gamma - \mu\gamma$.  Let $(1,\sigma_1^{'},\ldots, \sigma_n^{'},\epsilon^{'},\sigma_{n+1}^{'},\ldots,\sigma_{2n}^{'})$ be the evaluation dual basis for $K.(V) \otimes \Q$ and observe
that
\begin{equation} \nn
(1,\sigma_1^{'}, \ldots, \sigma_{n-1}^{'}, \sigma_n^{'} + \epsilon^{'}, -2\epsilon^{'},2\sigma_{n+1}^{'} - \epsilon^{'}, 2\sigma_{n+2}^{'}+\sigma_{n+1}^{'}, \ldots, 2\sigma_{2n}^{'} + \sigma_{2n-1}^{'})
\end{equation}
is a basis for $K.V$.  These basis elements (excepting the initial 1) map to the generators of $K.\Omega Spin(2n+2)$ listed above.  The twisting map $K.\Omega Spin(2n+2) \ra (K.\ast)_{\tau}$ takes these generators in turn to the integers
\begin{equation} \nn
\left(\binom{k}{1}, \ldots,\binom{k}{n-1},\binom{k}{n},0, 2\binom{k}{n+1},
2\binom{k}{n+2}+\binom{k}{n+1},\ldots, 2\binom{k}{2n}+\binom{k}{2n-1}
\right).
\end{equation}

As we did for the odd $Spin$ groups, we can simplify the presentation of the $K$-homology of
the loop spaces of these even $Spin$ groups before embarking on the $\Tor$ computation.  Abbreviate
the generators of $K.\Omega Spin(2n+2)$ as $(a_1,\ldots,a_{n-1}, \hat{a}_n, b,
\check{a}_{n+1}, a_{n+2},\ldots,a_{2n})$ respectively.  In $K.\Omega Spin(4n+2)$, the relation
$\rho_n$ eliminates the generator $b$, while $\rho_i$ eliminates $a_{2i}$ for $n<i<2n$, and the relation 
$\rho_{2n} - \epsilon^2$ eliminates $a_{4n}$.  This leaves
\begin{equation} \nn
K.\Omega Spin(4n+2) = \frac{\Z[a_1,\ldots, a_{2n-1},\hat{a}_{2n}, \check{a}_{2n+1}, a_{2n+3},\ldots,
a_{4n-1}]}{(\rho_1,\ldots,\rho_{n-1})}.
\end{equation}
Similarly in $K.\Omega Spin(4n)$, the relation $\rho_n$ eliminates the generator $\check{a}_{2n}$, while $\rho_i$ eliminates $a_{2i}$ for $n<i<2n-1$, and $\rho_{2n-1} - \epsilon^2$ eliminates $a_{4n-2}$, leading to the presentation
\begin{equation} \nn
K.\Omega Spin(4n) = \frac{\Z[a_1,\ldots,a_{2n-2},\hat{a}_{2n-1}, b, a_{2n+1}, a_{2n+3}, \ldots,
a_{4n-3}]}{(\rho_1,\ldots,\rho_{n-1})}.
\end{equation}
The crucial subring $R_n = \Z[a_1, \ldots, a_{2n-2}]/(\rho_1,\ldots,\rho_{n-1})$ that we considered
for the odd $Spin$ groups is precisely the relation subring of both $K.\Omega Spin(4n+2)$ and $K.\Omega Spin(4n)$.  As such our previous $\Tor$ calculation carries over without modification.  This completes the proof of Theorem~\ref{thmtkh} for the even $Spin$ groups and the resulting even
$Spin$ cyclic orders are recorded in Theorem~\ref{thmcyc}.

\subsection{Poincare-Dual Bases and the Cyclic Order of $K^{\tau} G$}
\label{cycgen}

We describe a general procedure for computing the cyclic order of
$K^{\tau} G$ for any simple $G$ and illustrate the method
with the group $G_2$.  The referee has pointed out that Braun~\cite{braun}
conjectures that this cyclic order is the greatest common divisor
of the representations generating the Verlinde ideal.  This conjecture is
almost certainly correct and therefore provides an efficient (as compared
with the method below) means of computing these cyclic orders.

Let $V=G/H$ denote the Bott generating variety for $\Omega G$;
recall that the $K$-cohomology of $V$ is $R[H]/i^*I[G]$ where
$i:H \ra G$ is the inclusion.  Pick a module basis $\{w_i\}$ for this
ring and consider the Poincare-dual basis $\{Dw_i\}$ of $K.V$.
The image of this basis in $K.\Omega G$, which we also denote by
$\{Dw_i\}$, is a set of algebra generators for $K.\Omega G$.  Note
that for any set $\{y_i\}$ of algebra generators for $K.\Omega G$,
the cyclic order of $K^{\tau} G$ is given by
$\gcd\{\tau_k(y_i)-\tau_0(y_i)\}$, where $\tau_k$ and $\tau_0$
are respectively the twisted and untwisted maps from $K.\Omega G$
to $K.*$.  In section~\ref{tihi} we saw that $\tau_k(Dw_i)=
\lb L^k \cup w_i,[V] \rb$ where $L$ denotes the generating line
bundle on $V$.  Decompose $w_i$ into a sum of irreducible
representations $\sum v_{ij}$, and let $h(v_{ij})$ denote the
highest weight corresponding to $v_{ij}$.  The product $L^k \cup v_{ij}$ 
is again irreducible, with highest weight $kL+h(v_{ij})$, and
theorem~\ref{botthi} therefore applies: the pushforward
$\lb L^k \cup v_{ij},[V] \rb$ is either $0$ or is (plus or minus) the
dimension of the irreducible representation of $G$ with highest weight
$T(kL+h(v_{ij})+\rho)-\rho$, where $T$ reflects a weight into
the fundamental Weyl chamber.  This procedure expresses the
cyclic order of $K^{\tau} G$ as the greatest common
divisor of a finite set of differences of dimensions of
irreducible representations of $G$.

Recall that the Bott generating variety for $G_2$ is $G_2/U(2)$
for the short-root inclusion of $U(2)$.  Let $a$ and $b$ denote the 
fundamental weights of $G_2$ corresponding to the 7 and 14 dimensional 
representations; in particular $R[G_2]=\Z[a,b]$.
Similarly $R[U(2)]=\Z[f,t,t^{-1}]$, where $f$ and $t$
are respectively the standard representation and the determinant
representation.  The restriction map is
\begin{align}
i^*(a) &= f + f^2 t^{-1} - 1 + f t^{-1} \nn \\
i^*(b) &= t + f^3 t^{-1} - 2f + f^2 t^{-1} + f^3 t^{-2} - 2ft^{-1} + t^{-1}.\nn
\end{align}
Let $s=t^{-1}$; the $K$-cohomology of the generating variety is then
$$K^{\bdot}V = \Z[f,s]/(f+f^2s-1+fs,1+f^3s^2-2fs+f^2s^2+f^3s^3-2fs^2+s^2).$$
(Note that the description of $K^{\bdot}V$ in Clarke~\cite{clarke} omits
certain relations, as the ring given there is not finitely generated.)
An integral basis for $K^{\bdot}V$ is then $\{1, s, s^2, f, fs, f^2\}$.
These representations of $U(2)$ are irreducible except for $f^2$ which
splits as $(f^2-t)+t$.  

\begin{figure}[t]
\begin{center}
\caption{The Weight Lattice of $G_2$} \vspace{6pt}
\epsfig{figure=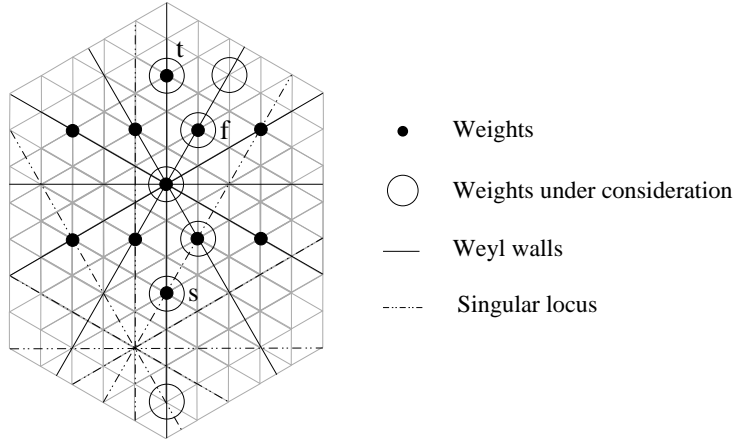,height=2.3in} \label{g2weights}
\end{center} \vspace{16pt}
\end{figure}

Consider the diagram of weights in Figure~\ref{g2weights}.
\nid The solid lines are the Weyl walls, the dotted lines describe the
set of singular weights, and the seven highest weights $h_i$ under 
consideration, namely $\{0,s,2s,f,f+s,2f,t\}$, are circled.  
Note that $L=t$ and as such, for $k>0$,
the weight $kL+h_i$ is either singular or is already in the fundamental
Weyl chamber.  The basis for $K.V$ is of course 
$\{D1,Ds,D(s^2),Df,D(fs),D(f^2)\}$ and we are interested in the
differences $\tau_k(Dw)-\tau_0(Dw)$.  Letting $\Gamma_{(n,m)}$ denote the 
dimension of the irreducible representation of $G_2$ with highest weight 
$na+mb$, the six differences are respectively
\vspace{-2pt}
\begin{align}
& \Gamma_{(0,k)} - \Gamma_{(0,0)} \nn\\
& \Gamma_{(0,k-1)} - 0 \nn\\
& \Gamma_{(0,k-2)} - 0 \nn\\
& \Gamma_{(1,k)} - \Gamma_{(1,0)} \nn\\
& \Gamma_{(1,k-1)} - 0 \nn\\
& \Gamma_{(2,k)} + \Gamma_{(0,k+1)} - \Gamma_{(2,0)} - \Gamma_{(0,1)}. \nn
\end{align} \newpage
Applying the Weyl dimension formula, we arrive at the cyclic order
\vspace{-2pt}
\begin{align}
c(G_2,k) = \gcd\{
& k(422+585k+400k^2+135k^3+18k^4)/120, \nn\\
& k(2+15k+40k^2+45k^3+18k^4)/120, \nn\\
& k(2-15k+40k^2-45k^3+18k^4)/120, \nn\\
& k(601+660k+350k^2+90k^3+9k^4)/30, \nn\\
& k(16+60k+80k^2+45k^3+9k^4)/30, \nn\\
& k(2867+2550k+1090k^2+225k^3+18k^4)/30
\}. \nn
\end{align}
Indeed, this agrees with the result from section~\ref{genSUnG2}.

\section{Twisted $\Spinc$ Bordism and the Twisted Index} \label{bord}

The ordinary $K$-homology of a space $X$ is entirely determined by
the $\Spinc$ bordism of $X$; see~\cite{hopkhov}.  This suggests that
much of the structure in twisted $K$-homology ought to be visible
in twisted $\Spinc$ bordism.  In section~\ref{tatesec} we saw that the 
cyclic order of the twisted $K$-homology of a group $G$ is determined 
by a collection of relations of the form $\tau_k(x)-\tau_0(x)=0$, where
$\tau_j$ is the $j$-twisted map from $K.\Omega G$ to $K.*$.  When the
class $x \in K.\Omega G$ is represented as the image of the fundamental
class of a $\Spinc$ manifold $M$, there is a natural $\Spinc$ manifold
$M(j)$ such that the fundamental class $[M(j)] \in MSpin^{c}_{\bdot}\ast$ 
maps via the
index to the element
$\tau_j(x) \in K.\ast$.  Moreover, there is an explicitly
identifiable twisted $\Spinc$ nullbordism over $G$ of $M(k)-M(0)$.  
In short, the
relations determining the cyclic order of twisted $K$-homology have
realizations in twisted $\Spinc$ bordism.  The construction of these
nullbordisms is the focus of sections~\ref{cocycle} and~\ref{nullbord}.
Section~\ref{exterior} discusses the possibility of representing the exterior 
generators of the twisted $K$-homology of $G$ by twisted $\Spinc$ manifolds.

\subsection{A Cocycle Model for Twisted $\Spinc$ Bordism} \label{cocycle}

In order to describe twisted $\Spinc$ manifolds explicitly, we need a
more geometric, less homotopy-theoretic description of twisted $\Spinc$ 
structures; in particular we present a cocycle model for twisted $\Spinc$
bordism.  This model is presumably well known and in any case takes cues 
from the Hopkins-Singer philosophy of differential functions~\cite{hopksing}.

Recall that $\Spinc$ is the total space of a $U(1)$-principal bundle over $SO$.
Correspondingly there is a principal bundle $BU(1) \ra B\Spinc \ra BSO$ which
is classified by $\beta w_2: BSO \ra BBU(1)$, the integral Bockstein of the 
second Stiefel-Whitney class.  A $\Spinc$ structure on an oriented manifold 
$M$ is a lift to $B\Spinc$ of the classifying map $\nu:M\ra BSO$ of the
(stable) normal bundle of $M$.  Such a lift is determined by a nullhomotopy 
of the composite $\beta w_2(\nu):M \ra BSO \ra BBU(1)$. 
Specifying such a nullhomotopy is equivalent to choosing a 2-cochain $c$ on 
$M$ such that the coboundary of $c$ is $\beta w_2(\nu(M))$.  (Note that
we have chosen 
once and for all a 3-cocycle $g$ representing the generator of 
$H^3(BBU(1);\Z)$, and the condition on the cochain $c$ is that $\delta c = 
\nu^*((\beta w_2)^*(g))$).  Ordinary $\Spinc$ bordism of $X$ is therefore
equivalent to bordism of oriented manifolds $M$ over $X$ equipped with a 
2-cochain $c$ on $M$ such that $$\delta c=\beta w_2(\nu(M)).\nn$$

The model for twisted $\Spinc$ bordism is similar.  We first recall the 
homotopy-theoretic definition of twisted $\Spinc$ bordism from 
section~\ref{tht}.  Given a twisting map $\tau: X \ra K(\Z,3)$, we
have a $K(\Z,2)$-principal bundle $P$ on $X$ and so an associated $B\Spinc$
bundle $Q=P \times_{K(\Z,2)} B\Spinc$.  More particularly we have a
series of bundles $$Q_n = P \times_{K(\Z,2)} B\Spinc(n)\nn$$ and
universal vector bundles $$UQ_n = (P \times_{K(\Z,2)} E\Spinc(n)) 
\times_{\Spinc(n)} \R^n.\nn$$  The corresponding Thom spectrum
$$\Thom(UQ) = P_+ \sm_{K(\Z,2)_+} M\Spinc \nn$$ has as its homotopy groups
the twisted $\Spinc$ bordism groups of $X$.  The twisted index map
to twisted $K$-homology is induced by the map $\id \sm \ind:
P_+ \sm_{K(\Z,2)_+} M\Spinc \ra P_+ \sm_{K(\Z,2)_+} K$.

The principal bundle $P$ and the associated $B\Spinc$ bundle $Q$ are defined 
by the pullbacks
\begin{equation} \nn
\xymatrix@C+5pt{
P \ar[r] \ar[d] & EK(\Z,2) \ar[d] & 
Q \ar[r] \ar[d] & EK(\Z,2) \times_{K(\Z,2)} B\Spinc \ar[d] \\
X \ar[r] & K(\Z,3) &
X \ar[r] & K(\Z,3).
}
\end{equation}
On the other hand $BSO$ is precisely the quotient $* \times_{K(\Z,2)}
B\Spinc$, and the diagram
\begin{equation} \nn
\xymatrix{
Q \ar[r] \ar[d] & BSO \ar[d]^{\beta w_2} \\
X \ar[r]_-{\tau} & K(\Z,3)
}
\end{equation}
is therefore a homotopy pullback.  Twisted $\Spinc$ bordism is the
homotopy of $\Thom(UQ)$; a map from a sphere into $\Thom(UQ)$ transverse
to the zero section $Q$ determines a manifold $M$ equipped with a
map $M \ra Q$.  This map $M \ra Q$ specifies
maps $i:M \ra X$ and $\nu:M \ra BSO$ (classifying
the normal bundle of $M$) and a chosen homotopy between $\tau i$ and 
$\beta w_2 \nu$.  The choice of this homotopy is equivalent to
the choice of a 2-cochain $c$ with coboundary equal to the difference
$\nu^*((\beta w_2)^*g) - i^*(\tau^*g)$, where $g$ is as before a 
3-cocycle representing the generator of the third cohomology of $K(\Z,3)$.
In summary, the $\tau$-twisted $\Spinc$ bordism of $X$ is bordism of oriented
manifolds $M$ equipped with a map $i:M \ra X$ and a 2-cochain $c$
such that $$\delta c = \beta w_2 (\nu(M)) - i^*(\tau),\nn$$ where 
$\nu(M)$ is the stable normal bundle of $M$.

\subsection{Twisted Nullbordism and the Geometry of the Cyclic Order} 
\label{nullbord}

In section~\ref{tatesec} we saw that the cyclic order of the
twisted $K$-homology of G is the greatest common divisor of the
collection of differences $\{\tau_k(x_i)-\tau_0(x_i)\}$, where 
$\{x_i\}$ is a set of algebra generators for $K.\Omega G$ and
$\tau_j$ denotes the $j$-twisted map from $K.\Omega G$ to $K.*$.
Frequently, these
generators $\{x_i\}$ can be described as the images of the fundamental
classes of $\Spinc$ manifolds $M_i$.  (For example, we gave such
a description for $SU(n+1)$, $Sp(n)$, and $G_2$ in sections~\ref{genSUnG2}
and~\ref{genSpn}).  In this case the manifolds $M_i$ admit modified
$\Spinc$ structures $M_i(j)$ and the index of 
$M_i(j)$ is $\tau_j(x_i) \in K.\ast$.  Moreover, there is a 
twisted $\Spinc$ structure (over $G$) on $M_i \times I$ cobounding the
difference $M_i(k) - M_i(0)$; the relations $\tau_k(x_i) - \tau_0(x_i) = 0$
determining the cyclic order of twisted $K$-homology therefore have 
realizations in twisted $\Spinc$ bordism.

Before constructing these twisted $\Spinc$ bordisms, we recall that
a $\Spinc$ structure can be altered by a line bundle and we discuss how
this alteration affects the pushforward of the fundamental class.
A twisted $\Spinc$ manifold is, as before, an oriented
manifold $M$ together with a 2-cochain $c$ such that $\delta c =
\beta w_2 (\nu(M)) - i^*(\tau)$.  In the examples we consider, the
underlying manifold $M$ is almost complex and so has a canonical ordinary
$\Spinc$ structure; in particular $M$ comes equipped with a 2-cochain
$b$ such that $\delta b = \beta w_2 (\nu(M))$.  A twisted structure on
$M$ is then given by a choice of 2-cochain $d$ such that $\delta d =
-i^*(\tau)$.  If the twisting class $\tau$ is zero on $M$, then the 
`twisted' $\Spinc$ structure corresponding to the cochain $b+d$ is
of course ordinary, but it nevertheless differs from the $\Spinc$
structure determined by the original cochain $b$.
We denote by $M(d)$ 
this modification of the canonical $\Spinc$ structure on $M$ by the 2-cocycle
$d$; we also refer to the alteration
as a modification by the corresponding line bundle
$L(d)$.  Let $\pi:M \ra *$ be the projection to a point; the pushforward
in $K$-theory depends on the $\Spinc$ structure on $M$ as follows:
\begin{equation}
\pi^{(M(d))}_!(1) = \pi^{(M)}_!(L(d)). \nn 
\end{equation}
This relation follows more or less immediately from the fact that the Thom 
class defined by the $\Spinc$ structure on $M(d)$ is $L(d)$ tensor the Thom 
class defined by the structure on $M$; see~\cite{spingeo}.  In terms of our
Rothenberg-Steenrod spectral sequence approach to twisted $K$-homology, this 
tells us that the twisted image $\tau_k([M])=\lb L(k),[M] \rb$ of the 
fundamental class $[M]$ is
reinterpretable as the ordinary image $\tau_0([M(k)])=\lb 1,[M(k)] \rb$ 
of the fundamental class $[M(k)]$.

\begin{center} $\sim\sim\sim\sim\sim$ \end{center} \vspace{6pt}

We present twisted $\Spinc$ bordisms realizing the relations in the
twisted $K$-homology of $SU(n+1)$.  In section~\ref{genSUnG2} we saw
that the fundamental classes $\{[\CP^i]\}$ are algebra generators 
for $K.\Omega SU(n+1)$, and the relations determining the cyclic order
of $K^{\tau} SU(n+1)$ are 
$$0 = \lb L^k,[\CP^i] \rb -1=\lb L^k,[\CP^i] \rb - \lb 1,[\CP^i] \rb . \nn$$
(These classes take values in the twisted $K$-theory of $SU(n+1)$ via
the inclusion $\CP^i \ra \Sigma \CP^i \ra \Sigma \CP^n \ra SU(n+1)$,
which is of course nullhomotopic.)
By the above remarks we can rewrite the relation as
$$0 = \lb 1,[\CP^i(k)]\rb -\lb 1,[\CP^i]\rb = \lb 1,[\CP^i(k) - \CP^i]\rb . 
\nn$$
If we can produce a nullbordism of $\CP^i(k) - \CP^i$, we
will have pulled the given relation back to twisted $\Spinc$ bordism.

Write $\Sigma \CP^i = ([-2,2] \times \CP^i) / (\{-2,2\} \times \CP^i \cup
[-2,2] \times *)$ and consider the inclusion $i:\Sigma\CP^i \ra \Sigma
\CP^n \ra SU(n+1)$.  Choose a 3-cocycle representing the twisting
$\tau$ on $SU(n+1)$ such that $i^*(\tau)$ is $k$ times the cocycle locally 
Poincare dual to the submanifold $\CP^{i-1} \times \{0\}$ in $\Sigma \CP^i$.  
The product $\CP^i \times [-1,1]$ has a canonical $\Spinc$ structure coming 
from the complex structure of $\CP^i$.  There is a twisted structure on 
$\CP^i \times [-1,1]$ defined by the 2-cochain 
$d$ that is $k$ times the cochain locally Poincare dual to the 
submanifold $\CP^{i-1} \times [-1,0]$; denote this twisted structure by
$(\CP^i \times [-1,1])(k|0)$.  The coboundary of $d$ is precisely 
$-i^*(\tau)$.  Moreover, the cochain $d$ restricts to $k$ times the 
generator
of $H^2(\CP^i \times \{-1\})$ and to zero on $\CP^i \times \{1\}$.  The
difference $\CP^i(k) - \CP^i$ is therefore zero in $MSpin^{c,\tau} SU(n+1)$,
as desired.  Notice that the same argument shows that \mbox{$\CP^i(l+k) - \CP^i(l)$}
is null for any $l$, which implies that $\binom{l+k+i}{i} - \binom{l+i}{i}$
is zero in $K^{\tau} SU(n+1)$.  In fact, for any sequence of integers
$\{l_i\}$, $1 \leq i \leq n$, the gcd of the set 
$\{\binom{l_i+k+i}{i}-\binom{l_i+i}{i}\}$ is again the cyclic order
of $K^{\tau} SU(n+1)$.

Whenever algebra generators of $K.\Omega G$ are represented as the fundamental
classes of $\Spinc$ manifolds, the same argument produces nullbordisms in
$MSpin^{c,\tau} G$ realizing the cyclic order of $K^{\tau} G$; we forgo
details.  Note though that in general the twisting cochain $d$ will no longer 
be locally Poincare dual to a submanifold but merely to an appropriate
singular chain.

\subsection{Representing the Exterior Generators of Twisted $K$-Homology} 
\label{exterior}

We would like to represent the algebra generators of $K^{\tau} G$
as the fundamental classes of twisted $\Spinc$ manifolds over $G$.  Here
we merely suggest an approach for further investigation, taking cues from 
the structure of the 
Rothenberg-Steenrod spectral sequence; in the process we produce a candidate 
representative for the generator of $K_1^{\tau} SU(3)$.  Finding
representatives in general will require a more thorough investigation of
$MSpin^{c,\tau}G$ and of the associated map to $K^{\tau}G$.

The structure of the ordinary $\Spinc$ bordism group is governed by
$\Spinc$ characteristic numbers; we briefly recall how to compute these
invariants.  In section~\ref{cocycle} we considered the principal bundle
$BU(1) \ra B\Spinc \ra BSO$ classified by $\beta w_2: BSO \ra BBU(1)$.
There is another principal bundle $B\Z/2 \ra B\Spinc \ra BSO \times BU(1)$
classified by $(w_2 \times r):BSO \times BU(1) \ra BB\Z/2$, where 
$r$ is the nontrivial map $BU(1) \ra BB\Z/2$.  This latter bundle is
usually more convenient for computations of $\Spinc$ characteristic classes.
The relationship between the two bundles is encoded in the matrix \pagebreak
\begin{equation} \nn
\xymatrix{
U(1) \ar[r]^{2} \ar[d] & U(1) \ar[r]^{r} \ar[d] & B\Z/2 \ar[r]^{\beta} \ar[d]
& BU(1) \ar[d] \\
\Spinc \ar[r] \ar[d] & \ast \ar[r] \ar[d] & B\Spinc \ar[r]^{\id} \ar[d]^{\pi
\times \lambda} & B\Spinc \ar[d]^{\pi} \\
SO \ar[r]^-{0} \ar[d]^{\Omega \beta w_2} & BU(1) \ar[r]^-{i} \ar[d]^{\id} 
& BSO \times BU(1) \ar[r]^-{\pi} \ar[d]^{w_2 \times r} 
& BSO \ar[d]^{\beta w_2} \\
BU(1) \ar[r]^{2} & BU(1) \ar[r]^{r} & BB\Z/2 \ar[r]^{\beta} & BBU(1).
}
\end{equation}
Indeed this diagram shows that the total spaces of the two fibrations
are the same.

Following Anderson, Brown, and Peterson~\cite{abp},
Stong~\cite{stong} showed that a $\Spinc$ manifold $M$ is zero in $\Spinc$ 
bordism if and only if all of its rational and mod 2 characteristic numbers 
vanish.  The map $$(\pi \times \lambda)^*:H^*(BSO \times BU(1);\Q) \ra
H^*(B\Spinc;\Q)\nn$$ is an isomorphism and $$(\pi \times \lambda)^*:H^*(BSO 
\times BU(1);\Z/2) \ra H^*(B\Spinc;\Z/2)\nn$$ is an epimorphism.  In
particular a $2n$-dimensional $\Spinc$ manifold $M$ is nullbordant if
all the characteristic classes of the underlying oriented manifold vanish
and the single $\Spinc$ characteristic number $\lb \lambda(M)^n,[M]_H\rb$ 
is zero. The characteristic class $\lambda$ depends on the $\Spinc$ structure 
on $M$ as follows.  Let $M(d)$ denote as in the last section the modification
of the $\Spinc$ structure on $M$ by the line bundle or 2-cocycle $d$.  The
class $\lambda(M(d))$ is then $\lambda(M)+2d$, as is easily checked by noting
that the composite $BU(1) \ra B\Spinc \overset{\lambda}{\ra} BU(1)$ is
multiplication by 2.

\begin{center} $\sim\sim\sim\sim\sim$ \end{center} \vspace{6pt}

We produce a candidate twisted $\Spinc$ representative for
the exterior generator of $K^{\tau} SU(3)$ by investigating
the corresponding class in the $E^2$ term of the Rothenberg-Steenrod
spectral sequence.  For simplicity we assume the twisting class $k$ is
odd; the even case is entirely analogous.

In section~\ref{torSUnSpn} we saw that the generator of $K_1^{\tau} SU(3)$
is represented at the $E^2$ term of the Rothenberg-Steenrod spectral
sequence by $x_2 - \frac{k+3}{2} x_1$; here $x_2$ and $x_1$ are elements
of $\Tor_1^{K.\Omega SU(3)}(\Z,\Z_{\tau})$, therefore of the $E^1$ term
of the spectral sequence, and their differentials are given by
\begin{align}
d^1 x_2 &= \lb 1,[\CP^2(k)]\rb - 1 = \lb 1,[\CP^2(k) - \CP^2(0)] \rb \nn \\
d^1 x_1 &= \lb 1,[\CP^1(k)]\rb - 1 = \lb 1,[\CP^1(k) - \CP^1(0)] \rb \nn
\end{align}
In section~\ref{nullbord} we found a twisted $\Spinc$ bordism
$X_2 = (\CP^2 \times I)(k|0)$ whose boundary has index
$$\ind(\partial X_2) = d^1 x_2.\nn$$
Because of this index property, we consider $X_2$ a geometric
representative of the algebraic class $x_2$.  Note that the bordism $X_2$ is
over $\Sigma \CP^2$ and therefore over $SU(3)$.

Similarly, we have a bordism $\widetilde{X_1} = (\CP^1 \times I)(k|0)$
whose boundary has index $\ind(\partial \widetilde{X_1}) = d^1 x_1$.
Given our selection of $X_2$, the manifold $\widetilde{X_1}$ is a
poor choice for a geometric representative of $x_1$; we
would like to have a five-dimensional bordism $X_1$, still living over
$\Sigma \CP^2$, with the same index property as $\widetilde{X_1}$.
A natural
choice for the underlying oriented bordism is $P(\nu + 1) \times I$,
where $P(\nu+1)$ is
the projectivization of the sum of a trivial bundle and the normal bundle 
of $\CP^1$ in $\CP^2$; this projectivization is a resolution of the Thom space
of the normal bundle and as such the bordism maps to $\CP^2 \times I 
\subset \Sigma\CP^2$.  There is moreover a twisted $\Spinc$ structure
on this bordism, denoted $X_1 = (P(\nu + 1) \times I)(k|0)$ and produced as 
in section~\ref{nullbord}, such that $$\ind(\partial X_1) = d^1 x_1.\nn$$

The linear combination $C=X_2 - \frac{k+3}{2} X_1$ wants to be an element of
$MSpin^{c,\tau} \Sigma \CP^2$ mapping to the exterior generator of
$K^{\tau} SU(3)$.  The trouble of course is that $C$ is not a closed
manifold and so does not properly represent an element of 
$MSpin^{c,\tau} \Sigma \CP^2$.  Note though that the map $\partial C
\ra \Sigma \CP^2$ is nullhomotopic by a nullhomotopy on which the
twisting class is zero.  Suppose there is a nullbordism $W$ of $\partial C$ 
in $M\Spinc *$; then the union $W \cup_{\partial C} C$ is a closed
twisted $\Spinc$ manifold over $\Sigma \CP^2$, as desired.  

The boundary of $C$ is $$\partial C = (\CP^2(k) - \CP^2) - \frac{k+3}{2}
(P(\nu+1)(k) - P(\nu+1)).\nn$$  All the $SO$-characteristic numbers 
of $\partial C$ certainly vanish.  The cohomology ring of $P(\nu+1)$ is
$H^{\bdot}(P(\nu+1)) = \Z[y,x]/(y^2,x^2+yx)$, where $y$ is the first Chern
class of the tautological bundle on the base $\CP^1$ and $x$ is the first 
Chern class of the fibrewise tautological bundle on the total space.
The tangential $\Spinc$ characteristic class of $P(\nu+1)(k)$ is 
$$ \lambda(T(P(\nu+1)(k))) = \lambda(T_{\text{horiz}}) + 
\lambda(T_{\text{vert}}) = -(2+2k)y - 2x, \nn$$
and the associated characteristic number is
$$ \lb \lambda(T(P(\nu+1)(k)))^2,[P(\nu+1)]_H \rb = 8k+4. \nn$$
Similarly the characteristic number for $\CP^2(k)$ is
$$ \lb \lambda(T(\CP^2(k)))^2,[\CP^2]_H \rb = 4k^2 + 12k + 9. \nn$$
The vanishing of the $\Spinc$ characteristic number for $\partial C$
follows:
$$ \lb \lambda(T(\partial C))^2,[\partial C]_H \rb = 4k^2 + 12k + 9 - 9
- \frac{k+3}{2} (8k+ 4 - 4) = 0. \nn$$
Picking any $\Spinc$ nullbordism $W$ of $\partial C$, the five-dimensional
twisted $\Spinc$ manifold $W \cup_{\partial C} C$ should represent the
generator of $K^{\tau}_1 SU(3)$.

\section*{Acknowledgments}

Marco Gualtieri prompted this inquiry with a question concerning the
twisted $K$-theory of $E_8$ and I thank him for the question and for
many fruitful conversations thereafter.  The algebraic geometry
interpretation of twisted $K$-theory mentioned in the introduction was
worked out in a series of entertaining discussions with Max 
Lieblich.  Francis Clarke kindly shared his notes on the $K$-homology
of the loop space of an even orthogonal group and thereby saved me
countless hours of frustration.  I would like to thank Nora Ganter, Teena 
Gerhardt, Andre Henriques, Mike Hill, Tyler Lawson, and Chris Woodward for helpful comments 
and for putting up with my various confusions along the way.  I also thank
the referee for writing such a careful report and for the valuable
remarks and corrections therein.  Mike Hopkins 
offered numerous indispensable suggestions and I am grateful for his 
guidance and for his ever astute and caring advice.

\end{document}